\documentclass[final]{siamltex}
\usepackage{xcolor}
\usepackage{amsmath}
\usepackage{amssymb}
\usepackage{mathrsfs}

\usepackage{graphicx}
 \usepackage{bm}
\usepackage{color}
\usepackage{float}
\usepackage{epstopdf}
\usepackage{caption}
\usepackage{subcaption}
\usepackage[breaklinks,colorlinks=true,linkcolor=blue,citecolor=red,backref=page]{hyperref}

\usepackage[textsize=small]{todonotes}

\DeclareMathAlphabet{\itbf}{OML}{cmm}{b}{it}
\DeclareMathAlphabet\mathbfcal{OMS}{cmsy}{b}{n}

\renewcommand{\tilde}{\widetilde}
\def\RR{\mathbb{R}}
\def\bx{{{\itbf x}}}

\def\by{{{\itbf y}}}

\def\bu{{{\itbf u}}}

\def\bg{{\itbf g}}
\def\be{{\boldsymbol{\epsilon}}}

\def\eps{{\epsilon}}

\def\bQ{{\itbf Q}}
\def\bv{{\itbf v}}
\def\bi{{\itbf i}}
\def\bet{{\boldsymbol{\eta}}}

\def\bV{{\itbf V}}
\def\bI{{\itbf I}}

\def\bF{{\itbf F}}
\def\bD{{\itbf D}}

\def\bff{{\itbf f}}
\def\cP{\mathcal{P}}
\def\cL{{\mathcal L}}

\def\bL{{\itbf L}}

\def\bA{{\itbf A}}

\def\cE{{\mathcal E}}

\def\bR{{\itbf R}}

\def\bpsi{{\boldsymbol{\psi}}}

\def\bphi{{\boldsymbol{\phi}}}
\def\bgamma{{\boldsymbol{\gamma}}}
\def\bPhi{{\boldsymbol{\Phi}}}
\def\bbM{\boldsymbol{\mathbb{M}}}
\def\bbS{\boldsymbol{\mathbb{S}}}

\def\FWI{{\scalebox{0.5}[0.4]{FWI}}}

\def\RM{{\scalebox{0.5}[0.4]{ROM}}}

\def\Gal{{\scalebox{0.5}[0.4]{GAL}}}

\def\cPR{{\boldsymbol{\cP}^{\RM}}}

\def\bphi{{\boldsymbol{\varphi}}}

\def\la{\lambda}
\def\12{{\frac{1}{2}}}

\DeclareMathAlphabet{\itbf}{OML}{cmm}{b}{it}
\DeclareMathAlphabet\mathbfcal{OMS}{cmsy}{b}{n}

\renewcommand{\tilde}{\widetilde}
\def\RR{\mathbb{R}}
\def\bx{{{\itbf x}}}

\def\by{{{\itbf y}}}

\def\bu{{{\itbf u}}}

\def\bg{{\itbf g}}
\def\be{{\boldsymbol{\epsilon}}}

\def\eps{{\epsilon}}

\def\bv{{\itbf v}}
\def\bi{{\itbf i}}
\def\bet{{\boldsymbol{\eta}}}

\def\bV{{\itbf V}}
\def\bI{{\itbf I}}

\def\bF{{\itbf F}}
\def\bD{{\itbf D}}
\def\bff{{\itbf f}}
\def\cP{\mathcal{P}}
\def\cL{{\mathcal L}}

\def\bL{{\itbf L}}

\def\bA{{\itbf A}}

\def\cE{{\mathcal E}}

\def\bR{{\itbf R}}

\def\bpsi{{\boldsymbol{\psi}}}

\def\bphi{{\boldsymbol{\phi}}}
\def\bPhi{{\boldsymbol{\Phi}}}
\def\bbM{\boldsymbol{\mathbb{M}}}
\def\bbS{\boldsymbol{\mathbb{S}}}

\def\FWI{{\scalebox{0.5}[0.4]{FWI}}}

\def\RM{{\scalebox{0.5}[0.4]{ROM}}}

\def\Gal{{\scalebox{0.5}[0.4]{GAL}}}

\def\cPR{{\boldsymbol{\cP}^{\RM}}}

\def\bphi{{\boldsymbol{\varphi}}}

\def\REG{{\scalebox{0.5}[0.4]{REG}}}

\def\la{\lambda}
\def\12{{\frac{1}{2}}}

\newtheorem{rem}[theorem]{Remark}
\newtheorem{algorithm}{Algorithm}
\newtheorem{assume}{Assumption}

\newcommand{\bc}{\textcolor{black}}

\graphicspath{{./FIGS/}}

\usepackage{showlabels}
\begin{document}

\title{Reduced order modeling for  \bc{first order} hyperbolic systems with application to multiparameter acoustic waveform inversion} \author{Liliana Borcea\footnotemark[1] \and Josselin
  Garnier\footnotemark[2] \and Alexander V. Mamonov\footnotemark[3] \and J\"{o}rn Zimmerling\footnotemark[4]}

\maketitle


\renewcommand{\thefootnote}{\fnsymbol{footnote}}
\footnotetext[1]{Applied Physics and Applied Mathematics, Columbia University, New York, NY,  10027. {\tt lb3539@columbia.edu}}
\footnotetext[2]{CMAP, CNRS, Ecole Polytechnique, Institut Polytechnique de Paris, 91120 Palaiseau, France.  {\tt
    josselin.garnier@polytechnique.edu}}
\footnotetext[3]{Department of Mathematics, University of Houston, TX 77204-3008. {\tt mamonov@math.uh.edu}}
\footnotetext[4]{Uppsala Universitet, Department of Information Technology, Division of Scientific Computing, 75105 Uppsala, Sweden. {\tt jorn.zimmerling@it.uu.se}}
\markboth{L. Borcea, J. Garnier, A.V. Mamonov, J. Zimmerling}{ROM forhyperbolic systems}

\begin{abstract}
Waveform inversion seeks to estimate  an inaccessible heterogeneous medium from data gathered by sensors that emit probing signals 
and measure the generated waves. It is an inverse problem for a second order wave equation or a first order hyperbolic system, with the sensor excitation modeled as a forcing term and the heterogeneous medium described by unknown, spatially variable coefficients. The traditional ``full waveform inversion" (FWI) formulation estimates the unknown coefficients via minimization of the nonlinear, least squares data fitting objective function.  For typical band-limited and high frequency data, this objective function has spurious local minima near and far from the 
true coefficients. Thus, FWI implemented with gradient based optimization algorithms may fail, even for good initial guesses. Recently, it was shown that it is possible 
to obtain a better behaved objective function for wave speed estimation, using data driven reduced order models (ROMs) that capture the propagation of pressure waves, governed by the classic second order wave equation. Here we introduce ROMs for vectorial waves,  satisfying a 
general first order hyperbolic system. They are defined via Galerkin projection on the space spanned by the wave snapshots, evaluated on  a uniform time grid with appropriately chosen time step. Our ROMs are data driven: They are computed in an efficient and non-iterative manner, from the sensor measurements, without knowledge of the medium and the snapshots. The ROM computation applies to any linear waves in lossless and non-dispersive media. For the inverse problem we focus attention on 
acoustic waves in a medium with unknown variable wave speed and density. We show that  these  can be determined via minimization of 
an objective function that uses a ROM based approximation of the vectorial wave field inside the inaccessible medium.  We assess the performance of our inversion approach  with numerical simulations and compare the results to those given by FWI.
\end{abstract}
\begin{keywords}
Hyperbolic systems, inverse wave scattering, data driven, reduced order modeling.
\end{keywords}

\begin{AMS}
65M32, 41A20
\end{AMS}

\section{Introduction}
\label{sect:intro}
Waveform inversion is an important technology in radar and sonar imaging \cite{curlander1991synthetic,cheney2009fundamentals,gilman2017transionospheric,blondel1997handbook}, seismology and geophysical exploration \cite{symes2008migration,virieux2009overview}, medical imaging with ultrasound \cite{szabo2004diagnostic}, and so on. It is an inverse problem concerned 
with the  estimation of an inaccessible,  heterogeneous medium, from time resolved measurements of the wave field, gathered with  source and receiver sensors.  

Mathematically, the heterogeneous medium is modeled by unknown variable coefficients in second order wave equations, or first order hyperbolic systems of equations. The forcing in these equations is localized at sources that emit probing signals, which are typically pulses.  The unknown coefficients depend on the type of waves: 
They are the bulk modulus and mass density for acoustic waves,  the dielectric permittivity and magnetic permeability for electromagnetic waves and  the Lam\'{e} parameters and density for elastic waves. 

Much of the existing literature is concerned with scalar (acoustic pressure) waves, modeled by the  second order  wave equation with constant density and unknown variable bulk modulus, which in turn determines the variable wave speed. This is a simplification, but the problem remains difficult because the mapping between the unknown wave speed and the pressure field at the receivers is complicated and nonlinear. In fact, for \bc{application relevant high frequency probing pulses}  and back-scattering data acquisition geometries, where the sources and receivers lie on the same side of the medium, the nonlinear, least squares data fitting objective function has numerous spurious minima, far and near the true wave speed.  This  behavior is known as ``cycle skipping". It is the main \bc{difficulty faced by the
``full waveform inversion" (FWI) approach implemented with gradient based minimization of the least-squares data misfit, even for reasonable initial guesses \cite{virieux2009overview}. } 

\bc{The  nonlinear least squares data fitting formulation of inversion has received a lot of attention since its introduction in geophysics  
in \cite{lailly1983seismic,tarantola1984inversion}. It was named  ``full waveform inversion" (FWI) 
in \cite{pan1988full}, to emphasize that it uses full seismogram information, not just arrival times. FWI is a computationally intensive, PDE constrained optimization 
that relies on application specific expertise in  algorithms, parametrization and regularization \cite{modrak2016seismic,barucq2019priori}.
Most studies consider the time domain formulation, but they may involve sequential frequency band filtering of the data, from low to high. 
Low frequencies help in estimating the kinematics (smooth part of the wave speed) 
which mitigates cycle skipping. There are also FWI methods that operate in the frequency domain, with low to high frequency sweeps and also with complex frequencies. 
The PhD thesis \cite{faucher2017contributions} is a comprehensive study of such methods. The research in FWI continues
and the methodology is used outside the geophysics community, for example in brain imaging with low frequency ultrasound \cite{cudeiro2022design}, ground penetrating radar \cite{lavoue20142d} and non-destructive evaluation of corroded materials \cite{rao2016guided}. 
 There remain challenges to overcome, especially for applications with high frequency waves, where the bandwidth is small with respect to the central frequency of the probing pulses. }

There are other waveform inversion approaches, \bc{beyond FWI. For example, the source type integral equation (STIE) method 
introduced in \cite{Habashy1994}  uses the observation that inhomogeneities in the medium introduce ``scattering source currents" that are related to the measured wave field by linear integral  operators.  Unfortunately, these operators have a large null space. 
To deal with the non-uniqueness, the STIE method uses multiple excitations and 
frequencies and it estimates the scattering currents and subsequently the material properties via a nonlinear iterative procedure. The STIE methodology has been extended to medical imaging with MRI in \cite{Remis2015}.}

\bc{Among other, more recent advances to waveform inversion, we mention the following three:} One approach
is to minimize the data misfit in the Wasserstein metric  \cite{EngquistFroese,yang2018application}. The resulting objective  function is improved; it is demonstrably convex for some simple models \cite{engquist2022optimal}, although the result does not hold for general media \cite[Fig. 3.3]{borcea2023data}. Another approach is the ``modeling operator extension" \cite{huang2018source} or ``extended FWI"
\cite{herrmann2013,warner2016}. It introduces additional degrees of freedom in the optimization, e.g., an unknown source signal in addition to the 
wave speed, and then systematically drives the optimization towards a meaningful result. This approach has been analyzed for  a transmission data acquisition geometry in a simple medium  in  \cite{symes2022error}, but guarantees for general media and back-scattering data
are lacking. A third approach,  introduced in \cite{druskin2016direct} for one dimensional media (see also the early study \cite{1099514}) and then developed further in \cite{druskin2018nonlinear,DtB,borcea2020reduced,borcea2021reduced,borcea2022internal,borcea2022waveform,borcea2023data} for more general media, computes from the data a reduced order model (ROM). This is a matrix  that captures essential features of wave propagation in the unknown medium and therefore, it can be used to  define objective functions that are amenable to gradient based minimization  \cite{borcea2022internal,borcea2022waveform,borcea2023data}. 

The advantage of FWI, extended FWI and the Wasserstein metric approaches is that they work with arbitrary placement of the sources and receivers. So far, the ROM methodology has been limited to full knowledge of  the ``array response matrix" gathered by a collection (active array) of coinciding sources and receivers. Active arrays are used for example in  radar and phased array ultrasonics, but are not feasible in all applications. Two other data acquisition setups, where the array response matrix can be obtained via some processing, are described in \cite[Section 4]{borcea2023data} and \cite{borcea2022waveform}. Nevertheless,  general source and receiver placements have not been addressed so far.  The advantage of the ROM approach, \bc{first documented in \cite{borcea2022waveform,borcea2023data}},  is that it can estimate strongly heterogeneous media for which the other approaches have failed.

In this paper we extend the ROM inversion methodology to  first order hyperbolic systems, that govern all types of linear waves in non-dispersive and lossless media. The ROM derivation is general, but we focus attention on multiparametric inversion for acoustic waves, where the wave speed $c$ and density $\rho$ are to be estimated from measurements of the vectorial wave field, with components given by the acoustic pressure and velocity. Like in \cite{druskin2018nonlinear,DtB,borcea2020reduced,borcea2021reduced,borcea2022internal,borcea2022waveform,borcea2023data}, the ROM is 
defined via Galerkin projection  \cite{benner2015survey,hesthaven2022reduced,brunton2019data} on the space spanned by the wave snapshots, sampled on a uniform time grid, with step $\Delta t$ chosen according to the Nyquist sampling criterion for the probing source signals.  
The snapshots  are vector valued fields that are unknown at points inside the medium. We show that they 
satisfy an exact time stepping scheme driven by a ``wave propagator" operator.  The ROM is a matrix of size determined by the number of source excitations and 
time steps up to the duration of measurements. It is  defined as the orthogonal projection of the wave propagator operator onto the unknown space of the snapshots. We derive the nonlinear mapping between the wave measurements at the sensors and the ROM and then obtain a 
data driven ROM construction via a non-iterative procedure that uses efficient tools from  numerical linear algebra.  

The results in 
\cite{druskin2018nonlinear,DtB,borcea2020reduced,borcea2021reduced,borcea2022internal,borcea2022waveform,borcea2023data} are specialized 
for the second order wave operator $\partial_t^2 - c^2 \Delta$, and are based on the positive definiteness and self-adjointness of  $-c^2 \Delta$, 
in the $L^2_{c^{-2}}$ inner product, weighted by $c^{-2}$. For first order hyperbolic systems the wave operator is $\partial_t + \cL$,  with skew-adjoint, first order partial differential operator $\cL$. \bc{This changes the technical aspects of the ROM construction, which is new and different from the previous studies. In the second order formulation \cite{druskin2018nonlinear,DtB,borcea2020reduced,borcea2021reduced,borcea2022internal,borcea2022waveform,borcea2023data}, the 
wave field is scalar valued (acoustic pressure) and the  propagator  is a self-adjoint operator. Here, the wave field is vectorial and it can be expressed mathematically as a time-dependent flow, with unitary 
propagator operator $e^{-\Delta t \cL}$ defined in terms of the unitary semigroup generated by $\cL$. The data driven computation of the ROM is novel and the ROM propagator has a different algebraic structure. It was block  tridiagonal and symmetric in \cite{druskin2018nonlinear,DtB,borcea2020reduced,borcea2021reduced,borcea2022internal,borcea2022waveform,borcea2023data}, whereas here we obtain a block upper Hessenberg ROM propagator,  that can be interpreted as a special, upwind  finite difference scheme for the  causal wave propagation. 
From the practical point of view, the first-order formulation considered in this paper opens the door  to addressing vector waves such as elastic and electromagnetic and to explore the use of different types of (polarized) sources and receivers.}

We use the ROM to obtain an approximation of the vectorial wave field at the inaccessible points inside the medium, aka the ``internal wave". This approximation is consistent with the measurements but the internal wave does not solve the hyperbolic system of equations unless the estimates of $c$ and $\rho$ are right. This motivates our multiparametric inversion approach which minimizes the discrepancy between the internal wave and the solution of the hyperbolic system, throughout the domain. This is different from the FWI approach, which minimizes the misfit of the wave at the sensor locations.  As demonstrated by numerical simulations, our objective function is better suited for minimization with gradient based algorithms, like Gauss-Newton, than that of FWI. 

\bc{As far as we know, no existing waveform inversion methodology comes with a mathematical guarantee of converge to the true solution i.e., the objective function is not convex. 
We cannot prove convexity either, but what we can show is that the ROM can resolve better the medium from the same data. 
There are many technical ``tricks" that one can use to improve the performance of FWI \cite{faucher2017contributions,modrak2016seismic}.
Our goal is to compare how the same optimization algorithm, with the same data, parametrization of the unknowns and regularization strategy,
performs on our objective function vs. that of FWI. We point the reader to  \cite{borcea2022waveform,borcea2023data} for an explicit  comparison of the FWI and a ROM based objective functions in a two dimensional search domain for the velocity $c$. The performance  of the ROM introduced in this paper is similar  to that 
in \cite{borcea2022waveform,borcea2023data} for the case of a medium with constant density $\rho$. Thus, those results stand and demonstrate, 
at least in those simple cases, the advantage of the ROM based inversion. Here we are interested in the more difficult problem, where we search for both $c$ and $\rho$. }

The paper is organized as follows: We begin in section \ref{sect:GenROM} with the derivation of the data driven ROM for a general first order 
hyperbolic system. In section \ref{sect:Acoustics} we specialize the result to acoustic waves in a medium with unknown wave speed $c$ and density $\rho$. The multiparametric inversion  is discussed in section \ref{sect:Inversion}. We also give there numerical simulation results. We end with a brief summary in section \ref{sect:Summary}.

\section{ROM for first order hyperbolic systems}
\label{sect:GenROM}
Since we deal with a time domain formulation, where the measurements are recorded for a duration $T$, we use the hyperbolicity of the problem and the finite wave speed to truncate the domain of wave propagation, which is typically the whole space $\RR^d$, to $\Omega \subset \RR^d$, with $d \in \{2,3\}$. We let  $\Omega$ be a bounded and simply connected domain, with $T$ dependent diameter, and with piecewise smooth boundary $\partial \Omega$, modeled with some appropriate homogeneous boundary conditions\footnote{Because  $\Omega$ is defined via truncation, the choice of boundary conditions does not affect the measurements up to time $T$.}. For example, for acoustic waves, 
$\partial \Omega$ may be sound hard or sound soft. For electromagnetic waves, we may assume a perfectly conducting $\partial \Omega$.

All  linear waves in lossless and non-dispersive media are governed by a  first order hyperbolic system of the form
\begin{equation}
\partial_t \bpsi_{\be}(t,\bx) + \cL \bpsi_{\be}(t,\bx) = 
\bc{\bff_{\be}(t,\bx)}, 
\qquad t \in \RR, ~ \bx \in \Omega.
\label{eq:LS1}
\end{equation}
The dimension $m$ of the vectorial wave field $\bpsi_{\be} \in \RR^m$ depends on $d$ and the type of waves. 
\bc{For example, in $d$-dimensional acoustics, the components of $\bpsi_{\be}$ are the pressure and velocity, so $m = d+1$. In three-dimensional electromagnetics ($d=3$), the wave consists of the electric and magnetic fields, so $m = 6$. In  two-dimensional electromagnetics (a three-dimensional case with medium  invariant in one direction)   the governing equations decouple into two polarizations with $m=3$ each.}

The variable $t$ in \eqref{eq:LS1} denotes time and $\cL$ is a first order, skew-adjoint  partial differential operator with respect to $\bx$, acting on $m$-dimensional  vectorial functions in $\bL^2(\Omega)$,  with homogeneous boundary conditions.  The  
vector  $\be = (\eps_1,\eps_2) \in \cE\subset \mathbb{N}^2$ indexes the wave excitation, represented by  
the  $m$-dimensional vector function \bc{$\bff_{\be}$ in 
$\bL^2( (-t_s,t_s)\times\Omega)$}, 
that models a source localized\footnote{Localized means that the support is contained in a ball centered at $\bx_{\eps_1}$, with radius that is small with respect to the wavelength.} at point $\bx_{\eps_1}$, 
with polarization indexed by $\eps_2$. \bc{In some applications, the forcing $\bff_{\be}$ may be a separable function of $t$ and $\bx$. For generality, we allow the probing signals to change from one source to another. 
These signals are assumed to be pulses,  supported in time in the interval $(-t_s,t_s)$, with modulation at central (carrier) frequency $\nu$.}
Prior to the excitation, there is no wave,
\begin{equation}
\bpsi_{\be}(t,\bx) = {\bf 0}, \qquad t < -t_s, ~ \bx \in \Omega.
\label{eq:LS2}
\end{equation}

The unknown coefficients that model the heterogeneous medium are in $\cL$. The waveform inversion problem is to determine them 
from the time-dependent array response matrix  $\bA$ with entries 
\begin{equation}
\bc{
\bA_{\be',\be}(t) = \int_{-t_s}^{t_s} dt' \int_{\Omega} d \bx \, \left[\bff_{\be'}(-t',\bx)\right]^T \bpsi_{\be}(t-t',\bx).
}
\label{eq:LS3}
\end{equation}
Here $T$ denotes the transpose, while $\be = (\eps_1,\eps_2)$ and $\be' = (\eps'_1,\eps'_2)$  are 2-dimensional indexes in $\cE\subset\mathbb{N}^2$. Since $\bff_{\be'}$
is supported near the location $\bx_{\eps_1'}$ of one of the sensors, the measurements are approximated by  $\bpsi_{\be}$ evaluated there. That we consider this wave field along the force $\bff_{\be'}$ is a requirement of our data driven ROM construction. We show in  the next section that for acoustic waves, equation \eqref{eq:LS3}  models measurements of the velocity. 

%
\subsection{Time stepping and the snapshots}
\label{sect:Tstep}
Denote by $n_{\cE} $ the  cardinality of the set $\cE\subset \mathbb{N}^2$  of excitation indexes.
To simplify notation, we gather the solutions of \eqref{eq:LS1}, for all the excitations, in the matrix valued wave field
\begin{equation}
\bpsi(t,\bx) = \left(\bpsi_{\be}(t,\bx)\right)_{\be \in \cE},
\label{eq:LS4}
\end{equation}
with $m$ rows and $n_{\cE}$ columns. Similarly, we gather all the forces in 
\begin{equation}
\bc{\bff(t,\bx)= \left(\bff_{\be}(t,\bx)\right)_{\be \in \cE} \in \RR^{m \times n_{\mathcal{E}}}.}
\label{eq:LS5}
\end{equation}

The solution of \eqref{eq:LS1} can be written as 
\begin{equation}
\bc{\bpsi(t,\bx) = \int_{-t_s}^t dt' \,  e^{-(t-t') \cL} \bff(t',\bx),}
\label{eq:LS6}
\end{equation}
where $e^{-t \cL}$ is the unitary semigroup generated by $\cL$. This operator is understood to 
act columnwise on $\bff$. Let us shift the time by some  interval\footnote{We show later that, since $e^{-t \cL}$ is unitary, the ROM construction is invariant with respect to  $\tau$.} $\tau  \ge t_s$, 
and define the new, $m$-dimensional field 
\begin{equation}
\bphi(t,\bx) = \bpsi(t+\tau,\bx),
\label{eq:LS7}
\end{equation}
with initial state 
\begin{equation}
\bc{
\bphi(0,\bx) = \bphi_0(\bx) =  \int_{-t_s}^{t_s} dt' \, e^{-(\tau-t') \cL } \bff(t',\bx). }
\label{eq:LS8}
\end{equation}
 Here we used that  $\bff$ is supported in $(-t_s,t_s)$.  Equation \eqref{eq:LS6} implies that $\bphi$ satisfies the homogeneous 
 equation 
 \begin{equation}
 \partial_t \bphi(t,\bx) + \cL \bphi(t,\bx) = {\bf 0}, \qquad t > 0, ~ \bx \in \Omega,
 \label{eq:LS9}
 \end{equation}
 and it takes the form of the time-dependent flow 
 \begin{equation}
 \bphi(t,\bx) = e^{-t \cL } \bphi_0(\bx).
 \label{eq:LS10}
 \end{equation}

The time stepping equation follows easily from \eqref{eq:LS10}: Let $t_j = j \Delta t$ be a uniform time grid, with $j \in \mathbb{N} $ and 
$\Delta t\lesssim 1/(2 \nu)$, so that the Nyquist sampling criterion for the signal $s$ is satisfied. Define the wave snapshots  on this grid by
\begin{equation}
\bphi_j(\bx) = \bphi(t_j,\bx), \qquad j \ge 0.
\label{eq:LS11} 
\end{equation}
These satisfy  the time stepping iteration
\begin{equation}
\bphi_{j+1}(\bx) = \cP \bphi_j(\bx),  \qquad j \ge 0, ~ \bx \in \Omega, 
\label{eq:LS12}
\end{equation} 
driven by the wave propagator operator 
\begin{equation}
\cP = e^{-\Delta t \cL}.
\label{eq:LS13}
\end{equation} 
This unitary operator contains all the information about the medium and  the ROM approximates it by the matrix  $\cPR$, called the ``ROM propagator".

\subsection{Galerkin approximation of the snapshots}
\label{sect:Gal}
We are interested in the Galerkin projection of the time stepping equation  \eqref{eq:LS12} on the space $\mathcal{S}^{\Gal}$ spanned by the first $n_t$ 
snapshots, gathered in  the matrix valued field
\begin{equation}
\bPhi(\bx) = \left(\bphi_0(\bx), \ldots, \bphi_{n_t-1}(\bx) \right),
\label{eq:LS14}
\end{equation}
with $m$ rows and $n_t n_{\cE}$ columns. The ROM is defined using: 

\vspace{0.05in}
\begin{assume}
\label{as:1}
The columns of $\bPhi(\bx)$ are  linearly independent vector fields.
\end{assume}

\bc{
\begin{rem}
Assumption \ref{as:1} holds in general i.e., its converse is unlikely to occur. Indeed, suppose that the columns of $\bPhi$ were linearly dependent. Then,  there exists a linear combination $\bPhi\bgamma $ of these columns  that is vanishing in $\Omega$,
where  the coefficients are stored in the column vector $\bgamma  =(\gamma_{j,\be})_{j=0,\ldots,n_t-1, \be \in {\cal E}} $ of size $n_t n_{\cE}$.
Equivalently, it is possible to build the source 
\[\bff^{(\bg)}(t,\bx) = \sum_{j=1}^{n_t} \sum_{\be \in {\cal E}} \gamma_{j,\be}  \bff_\be(t-t_j - \tau,\bx),\]
supported in space on the sensor locations and in time in $(-(n_t-1) \Delta t-\tau-t_s,0)$ such that the resulting wave vector field vanishes in $\Omega$ at time $0$, and hence at any positive time.
 It is very unlikely that a null controllability problem can be solved by "chance" that way.
 \end{rem}}

\bc{We will see later that the Gramian matrix of the snapshots stored in $\bPhi$ plays a key role in the data driven computation of the ROM.  Assumption \ref{as:1} implies that the Gramian is nonsingular. This is sufficient for the theory in this section. However, to get a stable ROM computation, the Gramian must also be well-conditioned. In our experience, this occurs for arrays with sensor separation of about half the wavelength, calculated at frequency $\nu$, in the reference medium near the sensors, and with time sampling step $\Delta t \lesssim 1/(2\nu)$. If the Gramian is ill-conditioned or even singular (i.e., Assumption \ref{as:1} does not hold),
then  the ROM construction requires regularization. We discuss this in section \ref{sect:Inversion}. }

\vspace{0.05in}
%

The Galerkin approximation of the snapshots is 
\begin{equation}
\bphi_j^\Gal(\bx) = \bPhi(\bx) \bg_j \in \mathcal{S}^\Gal, \qquad j \ge 0, ~ \bx \in \Omega,
\label{eq:LS15}
\end{equation}
where $\bg_j \in \mathbb{R}^{n_t n_{\cE} \times n_{\cE}}$ are the matrices of Galerkin coefficients, calculated 
so that, when substituting \eqref{eq:LS15} into \eqref{eq:LS12}, we get a residual that is orthogonal to the approximation space.
We write  this orthogonality relation using a particular orthonormal basis of $\mathcal{S}^\Gal$, given by the Gram-Schmidt orthogonalization procedure,
summarized in the equation
\begin{equation}
\bPhi(\bx) = \bV(\bx) \bR, \qquad \int_{\Omega} d \bx \, \bV^T(\bx) \bV(\bx) = \bI_{n_t n_\cE},
\label{eq:LS16}
\end{equation}
where  $\bI_{q}$ denotes the $q \times q$ identity, for any $q \in \mathbb{N}$.  The Gram-Schmidt procedure is  good for our context, because it gives a causal basis. That is to say, if we write $\bV$, similar to \eqref{eq:LS14},
\begin{equation}
\bV(\bx) = \left(\bv_0(\bx), \ldots, \bv_{n_t-1}(\bx) \right),
\label{eq:LS17}
\end{equation}
 we have 
\begin{equation}
\bv_j(\bx) \in \mbox{span}\{\bphi_l(\bx), 0 \le l \le j\}, \qquad j = 0, \ldots, n_t-1,
\label{eq:LS19}
\end{equation}
with 
\begin{equation}
\int_{\Omega} d \bx \, \bv_j^T(\bx) \bv_l(\bx) = \delta_{j,l} \bI_{n_\cE}.
\label{eq:LS18}
\end{equation}
The causality \eqref{eq:LS19} is reflected in equation \eqref{eq:LS16} by the block upper triangular algebraic structure of the matrix $\bR \in \mathbb{R}^{n_t n_{\cE} \times n_t n_{\cE}}$, with 
$n_{\cE} \times n_{\cE}$ blocks. 

The equation for the Galerkin coefficients is 
\begin{equation}
\int_{\Omega} d \bx \, \bV^T(\bx) [\bPhi(\bx) \bg_{j+1} - \cP \bPhi(\bx) \bg_j ] ={\bf 0}, \qquad j \ge 0.
\label{eq:LS20}
\end{equation}
or, equivalently, 
\begin{equation}
 \bR^{-T} \bbM \bg_{j+1}  = \bR^{-T} \bbS \bg_j, \qquad j \ge 0,
\label{eq:LS21}
\end{equation}
with $\bR^{-T}$  the transpose of the inverse of $\bR$. 
Here we solved for $\bV$ in equation \eqref{eq:LS16} and introduced the 
Gramian, aka ``mass"   matrix 
\begin{equation}
\bbM = \int_{\Omega} d \bx \, \bPhi^T(\bx) \bPhi(\bx) \in \RR^{n_t n_\cE \times n_t n_\cE},
\label{eq:LS22}
\end{equation}
and the ``stiffness" matrix 
\begin{equation}
\bbS = \int_{\Omega} d \bx \, \bPhi^T(\bx) \cP \bPhi(\bx) \in \RR^{n_t n_\cE \times n_t n_\cE}.
\label{eq:LS23}
\end{equation}

\vspace{0.05in}
There are a few important observations for the Galerkin version  \eqref{eq:LS21} of the time stepping  equation: 

\vspace{0.05in}\begin{enumerate}
\itemsep 0.05in
\item The matrices $\bR$ and $\bbM$ are invertible due to Assumption \ref{as:1}, so all the Galerkin coefficients $\bg_j$ are uniquely defined by 
equation \eqref{eq:LS21},  in terms of $\bg_0$.  We choose 
\begin{equation}
\bg_0 = \bi_0,
\label{eq:LS23p}
\end{equation}
where $\{\bi_j\}_{j=0}^{n_t-1}$ denote the $n_t n_{\cE} \times n_{\cE}$ ``column blocks" of the   identity matrix
$\bI_{n_t n_{\cE}} = \left(\bi_0 , \ldots, \bi_{n_t-1} \right)$. The choice \eqref{eq:LS23p} ensures that the Galerkin approximation \eqref{eq:LS15}  satisfies exactly the initial condition. 
\item  The first $n_t$ matrices of Galerkin coefficients are trivial, meaning that  
\begin{equation}
\bg_j = \bi_j, \qquad 0 \le j \le n_t-1,
\label{eq:LS24}
\end{equation}
and therefore, 
\begin{equation}
\bI_{n_t n_{\cE}} = \left(\bg_0, \ldots, \bg_{n_t-1} \right) = \left(\bi_0, \ldots, \bi_{n_t-1} \right).
\label{eq:LS25}
\end{equation}
\item The block upper triangular matrix  $\bR$ in the Gram-Schmidt orthogonalization is the block Cholesky square root of the mass matrix 
\begin{equation}
\bbM = \bR^T \bR.
\label{eq:LS26}
\end{equation}
\end{enumerate}
Note that  equation \eqref{eq:LS24} is equivalent to saying that the Galerkin approximation \eqref{eq:LS15} is exact not only at  $t = t_0$, 
but also at  $t = t_j$, for $j = 1, \ldots, n_t-1$. This follows from the definition of the approximation space $\mathcal{S}^\Gal$ and the fact that \eqref{eq:LS12} implies 
\begin{equation}
\bphi_{j+1}(\bx) - \cP \bphi_j(\bx) = \bPhi(\bx) \bi_{j+1} - \cP \bPhi(\bx) \bi_j = {\bf 0}, \qquad 0 \le j \le n_t-2.
\label{eq:LS27}
\end{equation}
Equation \eqref{eq:LS26} is derived from the definition \eqref{eq:LS22} of the mass matrix and the Gram-Schmidt equation \eqref{eq:LS16},
\begin{equation}
\bbM = \int_{\Omega} d \bx \, \left[ \bV(\bx) \bR\right]^T \bV(\bx) \bR = \bR^T \int_{\Omega} d \bx \, 
\bV^T(\bx) \bV(\bx) \, \bR = \bR^T \bR.
\label{eq:LS28}
\end{equation}

\subsection{Galerkin projection ROM}

Equation \eqref{eq:LS21} is the algebraic, aka  ROM equivalent of \eqref{eq:LS12}. The ROM snapshots are the $n_t n_{\cE} \times n_{\cE}$ matrices defined by 
\begin{equation}
\bphi_j^{\RM} = \bR \bg_j, \qquad j \ge 0,
\label{eq:LS29}
\end{equation}
and due to the block Cholesky factorization \eqref{eq:LS26}, they satisfy the time stepping equation
\begin{equation}
\bphi_{j+1}^\RM = \cPR \bphi_{j}^\RM, \qquad j \ge 0,
\label{eq:LS30}
\end{equation}
driven by the ROM propagator matrix
\begin{equation}
\cPR = \bR^{-T} \bbS \bR^{-1}.
\label{eq:LS31}
\end{equation}
This is a projection ROM: According to \eqref{eq:LS15}-\eqref{eq:LS16}, the ROM snapshots are the projection 
 of the Galerkin approximation of the snapshots of the wave field, 
\begin{equation}
\bphi_j^{\RM} = \bR \bg_j = \int_{\Omega} d \bx \, \bV^T(\bx) \bPhi(\bx) \bg_j = \int_{\Omega} d \bx \, \bV^T(\bx) \bphi_j^\Gal(\bx), \qquad j \ge 0,
\label{eq:LS32}
\end{equation}
and from equation \eqref{eq:LS23} we get
\begin{equation}
\cPR = \int_{\Omega} d \bx \, \bR^{-T} [\bV(\bx) \bR]^T \cP  [\bV(\bx) \bR] \bR^{-1} = \int_{\Omega} d \bx \, \bV^T(\bx)  \cP  \bV(\bx).
 \label{eq:LS33}
\end{equation}

\subsection{Data driven ROM}
\label{sect:dataROM}
From the ROM definition above, it appears that its computation requires the matrix $\bPhi$ of snapshots. The next theorem states that the ROM is data driven. 

\vspace{0.05in} \begin{theorem}
\label{thm.1} 
The ROM can be computed directly from the array response matrix~\eqref{eq:LS3}, that defines the ``data matrices"
\begin{equation}
\bc{ \bD_j =   \bA(t_j)+\bA^T(-t_j)  , \qquad j \ge 0.}
\label{eq:LS35}
\end{equation}
The data matrices satisfy
\begin{align}
\bD_j &= \int_{\Omega} d \bx \, \bphi_0^T(\bx) \bphi_j(\bx)\label{eq:LS34_p}\\
 &= \big(\bphi_0^{\RM}\big)^T \bphi_j^\RM, \qquad j = 0, \ldots, n_t-1,
\label{eq:LS34}
\end{align}
where the second equality means that the ROM interpolates the measurements on the uniform time grid.
{The ROM  is computed from $\{\bD_j\}_{j=0}^{n_t}$ and its  propagator matrix \eqref{eq:LS31} has an unreduced, upper Hessenberg block structure.}
\end{theorem}

\vspace{0.05in}
\begin{proof}
The proof of  \eqref{eq:LS34_p} follows,  basically,  from energy conservation.  Recall definition \eqref{eq:LS7} of $\bphi$ and 
define the matrix 
\begin{equation}
\bQ(t;t_j) = \int_{\Omega} d \bx \, \bpsi^T(t+\tau,\bx) \bpsi(t+t_j+\tau,\bx).
\label{eq:Pf1}
\end{equation} 
The trace of $\bQ(\cdot; t_0 = 0)$ is  the total energy of the system at time $t + \tau$, for all excitations. According to  \eqref{eq:LS7}, we have 
\begin{equation}
\bQ(0;t_j) = \int_{\Omega} d \bx \, \bphi_0^T(\bx) \bphi_j(\bx).
\label{eq:Pf2}
\end{equation} 
The time derivative of \eqref{eq:Pf1}  is 
\begin{align} 
\partial_t \bQ(t;t_j) &=  \int_{\Omega} d \bx \, \left[ \partial_t \bpsi^T(t+\tau) \bpsi(t+t_j+\tau,\bx) + \bpsi^T(t+\tau) \partial_t \bpsi(t+t_j+\tau,\bx)\right] \nonumber \\
&
\bc{
\stackrel{\eqref{eq:LS1}}{=} \int_{\Omega} d \bx \, \left\{ \left[-\cL \bpsi(t+\tau,\bx)+  \bff(t+\tau,\bx) \right]^T \bpsi(t+t_j+\tau,\bx)  \right. \nonumber 
}
\\
&
 \qquad \quad 
\bc{
 \left. + 
\bpsi^T(t+\tau,\bx) \left[-\cL \bpsi(t+t_j + \tau,\bx) + \bff(t+t_j + \tau,\bx) \right] \right\} \nonumber 
}
\\
& 
\bc{
{=} \int_{\Omega} d \bx \, \left\{  \bff^T(t+\tau,\bx)\bpsi(t+t_j+\tau,\bx)   + 
\bpsi^T(t+\tau,\bx)   \bff(t+t_j + \tau,\bx)  \right\} ,
}
\label{eq:Pf3}
\end{align}
where we used that $\cL$ is skew-adjoint. Integrating in $t$ and evaluating at $t = 0$ we get
\begin{align} 
\bQ(0;t_j) &
\bc{= \int_{-\infty}^0 d t \int_{\Omega} d \bx \left\{  \bff^T(t+\tau,\bx)\bpsi(t+t_j+\tau,\bx)   + 
\bpsi^T(t+\tau,\bx)   \bff(t+t_j + \tau,\bx)  \right\}  \nonumber
} \\
\nonumber
& \bc{
= \int_{-\tau}^{\infty} \int_{\Omega}\bff^T(-t',\bx)\bpsi(t_j-t',\bx)   d \bx d t' }
\\
& 
\bc{\quad  + 
\int_{-\tau-t_j}^{\infty}  \int_{\Omega} \bpsi^T(-t_j-t'',\bx)   \bff(-t'',\bx)  d\bx dt'' ,} 
\label{eq:Pf4}
\end{align}
where we made the change of variables $t+\tau = -t'$ and $t+t_j+\tau = -t''$. Since the source is supported in time in the interval 
$(-t_s,t_s)$ and $\tau \ge t_s$, we observe from definition \eqref{eq:LS3} that the right hand side in \eqref{eq:Pf4} equals the matrix \bc{$\bD_j=   \bA(t_j )  + \bA^T(-t_j)$} defined in equation \eqref{eq:LS35}. The proof of \eqref{eq:LS34_p} follows from \eqref{eq:Pf2}.
%
%
%

Next, we show that the mass matrix is data driven. Writing its definition \eqref{eq:LS22} block by block we get, 
for $0 \le j \le l \le n_t-1$, that  
\begin{align}
\bbM_{j,l} &= \int_{\Omega} d \bx \, \bphi_j^T(\bx) \bphi_l(\bx) \nonumber \\
&\hspace{-0.1in}\stackrel{\eqref{eq:LS10}}{=} \int_{\Omega} d \bx \, \left[e^{-t_j \cL}\bphi_0(\bx)\right]^T e^{- t_l \cL} \bphi_0(\bx) \nonumber \\
&=\int_{\Omega} d \bx \, \bphi_0^T(\bx) e^{-(t_l-t_j) \cL} \bphi_0(\bx) \nonumber  \\
&\stackrel{\eqref{eq:LS10}}{=} \int_{\Omega} d \bx \, \bphi_0^T(\bx) \bphi_{l-j}(\bx)  \nonumber \\
&\stackrel{\eqref{eq:LS34_p}}{=} \bD_{l-j}. \label{eq:LS36}
\end{align}
Here we used that $e^{-t \cL}$ is a unitary semigroup, since $\cL$ is skew-adjoint.
It is obvious from the definition \eqref{eq:LS22} that $\bbM$ is symmetric, so the  blocks below the diagonal are 
\begin{equation}
\bbM_{j,l} = (\bbM_{l,j})^T, \qquad 0 \le l < j \le n_t-1.
\label{eq:LS37}
\end{equation}

The stiffness matrix is also data driven: Starting from definition 
\eqref{eq:LS23}, we get 
\begin{align}
\bbS_{j,l} &= \int_{\Omega} d \bx \, \bphi_j^T(\bx) \cP \bphi_l(\bx) \nonumber \\
&\hspace{-0.1in}\stackrel{\eqref{eq:LS12}}{=} \int_{\Omega} d \bx \, 
 \bphi_j^T(\bx) \bphi_{l+1}(\bx)  \nonumber \\
 &= \bD_{l+1-j},
 \label{eq:LS38}
 \end{align}
for $0 \le j \le l \le n_t-1$. Unlike the mass matrix, $\bbS$ is not symmetric. However, we can 
calculate its remaining blocks, indexed by  $0 \le l < j-1$, with $j = 2, \ldots, n_t-1$, as follows
\begin{align}
\bbS_{j,l} &= \int_{\Omega} d \bx \, \bphi_j^T(\bx) \cP \bphi_l(\bx) \nonumber \\
&= \int_{\Omega} d \bx \, [\cP^T \bphi_j(\bx)]^T \bphi_l(\bx) \nonumber \\
&= \int_{\Omega} d \bx \, 
 \bphi_{j-1}^T(\bx) \bphi_{l}(\bx)  \nonumber \\
 &= \left[\int_{\Omega} d \bx \, 
 \bphi_{l}^T(\bx) \bphi_{j-1}(\bx)\right]^T  \nonumber \\
 &= \bD^T_{j-1-l} .
 \label{eq:LS38p}
 \end{align}

To prove the data interpolation relations \eqref{eq:LS34}, we deduce from equations \eqref{eq:LS24} and \eqref{eq:LS29} 
that 
\begin{equation}
\bR = \left(\bphi_0^\RM, \ldots, \bphi_{n_t-1}^\RM \right).
\label{eq:LS40}
\end{equation}
Since $\bR$ is the block-Cholesky square root of $\bbM$, we have 
\begin{equation}
\bD_j \stackrel{\eqref{eq:LS36}}{=} \bbM_{0,j} = \bi_{0}^T \bbM \bi_j = \bi_{0}^T \bR^T \bR \bi_j = (\bphi_0^{{\RM}})^T \bphi_j^\RM,
\label{eq:LS41}
\end{equation}
for $j = 0,\ldots, n_t-1$, as stated in \eqref{eq:LS34}.

It remains to prove the algebraic structure of the propagator ROM: We use equation \eqref{eq:LS40} and the  iteration  \eqref{eq:LS30}:  If we denote the $n_{\cE} \times n_{\cE}$, no-zero blocks of $\bR$  by $\bR_{j,l}$, for  $0 \le j \le l \le n_t-1$,
we get from equation \eqref{eq:LS30} evaluated at $j = 0$ that the first column of blocks of $\cPR$ satisfies
\begin{equation*}
\bphi_1^\RM = \begin{pmatrix} \bR_{0,1} \\ \bR_{1,1} \\ {\bf 0} \\ \vdots \\ {\bf 0} 
\end{pmatrix} = \cPR \bphi_0^\RM = \begin{pmatrix} \cP^\RM_{0,0} \\ \cP^\RM_{1,0} \\ \cP^\RM_{2,0} \\ \vdots \\ \cP^\RM_{n_t-1,0} 
\end{pmatrix} \bR_{0,0}.
\end{equation*} 
Assumption \ref{as:1} ensures that  $\bbM$ and therefore $\bR$ are invertible, which makes all the diagonal blocks of $\bR$  invertible. 
Thus, we have 
\begin{equation*}
\cP^\RM_{l,0} \bR_{0,0} = {\bf 0} \quad  \leadsto \quad \cP^\RM_{l,0} = {\bf 0}, \qquad  2 \le l \le n_{t-1},
\end{equation*}
and 
\begin{equation*}
\cP^\RM_{1,0} \bR_{0,0} = \bR_{1,1} \quad  \leadsto \quad \cP^\RM_{1,0} \ne {\bf 0}.
\end{equation*}
Now proceed inductively for $j \ge 1$, under the hypothesis that for $j = 0, \ldots, J$, with $J \le n_t-2$, we have 
$\cP^\RM_{j+1,j} \ne {\bf 0}$ and 
\begin{equation}
\cP^\RM_{l,j} = {\bf 0}, \qquad j+2 \le l \le n_{t-1}.
\label{eq:LS43}
\end{equation}
Then, equation \eqref{eq:LS30} evaluated at $j = J$ gives  
\begin{equation*}
\bphi_{J+1}^\RM = \begin{pmatrix} \bR_{0,J+1} \\ \bR_{1,J+1} \\ \vdots \\ \bR_{J+1,J+1}  \\ {\bf 0} \\ \vdots \\ {\bf 0} 
\end{pmatrix}  = \begin{pmatrix}\cP^\RM_{0,J} \\ \cP^\RM_{1,J} \\ \cP^\RM_{2,J} \\ \vdots \\ \cP^\RM_{n_t-1,J} 
\end{pmatrix} \bR_{J,J}
+ 
\sum_{l=0}^{J-1} \begin{pmatrix}\cP^\RM_{0,l} \\ \cP^\RM_{1,l} \\ \cP^\RM_{2,l} \\ \vdots \\ \cP^\RM_{l+1,l}\\ {\bf 0} \\ 
\vdots \\ {\bf 0} \end{pmatrix} \bR_{l,J}.
\end{equation*}
Equating the entries in the $l^{\rm th}$ rows, with $l > J+1$, and using that $\bR_{J,J}$ is invertible, we deduce that the induction hypothesis equation \eqref{eq:LS43} extends to $J$, as well.  
We also have that 
\begin{equation*}
\cP^\RM_{J+1,J} \bR_{J,J} = \bR_{J+1,J+1} \quad  \leadsto \quad \cP^\RM_{J+1,J} \ne {\bf 0}.
\end{equation*}
The induction argument terminates at  $j = n_t-2$, and proves that  $\cPR$ has unreduced block upper Hessenberg structure. 
\end{proof}

\vspace{0.05in}\begin{rem} 
\label{rem:rec}
Definition \eqref{eq:LS35} shows that if there is reciprocity in the response of the medium i.e., if $\bA$ is symmetric, then 
the data matrices $\bD_j$ are symmetric and we can drop the transpose in equation \eqref{eq:LS38p}. We show in 
section \ref{sect:Acoustics} that  reciprocity holds 
in the case of acoustic waves.
\end{rem}

\vspace{0.05in}\begin{rem} 
\label{rem:rec1}
The second term in definition \eqref{eq:LS35} contributes only at the very early time instants $t_j \le 2 t_s$, due to causality.
Indeed, at later instants, satisfying $t_j > 2 t_s$, we have  
\bc{$\bpsi^T(-t_j-t,\bx) = {\bf 0}$}
 in the integral, because  $-t_j-t \in (-t_j-t_s,-t_j+t_s)$ does not intersect the time interval $(-t_s,\infty)$ where the wave is supported.
\end{rem}

\subsection{Properties of the ROM}
\label{sect:PropROM}
The calculations above reveal two important aspects of the ROM, that account for the action of the unitary semigroup $e^{-t \cL}$ generated by  $\cL$:

\vspace{0.05in} 1.  The mass matrix $\bbM$ is symmetric, with block Toeplitz structure. Its main diagonal contains the initial data matrix $\bD_0$, 
which is symmetric by equation \eqref{eq:LS34_p} evaluated at $j = 0$. The other diagonals contain the data matrices
 $\bD_j$, for $j = 1, \ldots, n_t-1$, respectively. The stiffness matrix $\bbS$ is also 
block Toeplitz, but its diagonals are shifted. The initial data blocks $\bD_0$ appear in the first sub-diagonal, the main diagonal 
contains the data matrix $\bD_1$, and the remaining upper diagonals have the blocks $\bD_j$, for $j = 2, \ldots, n_t$.

\vspace{0.05in} 2. The ROM propagator $\cPR$ preserves the unitary structure of the true wave propagator operator $\cP = \exp( -\Delta t \cL)$. 
To see this, let us write the iteration \eqref{eq:LS12} backward in time, using that $\cP$ is unitary, 
\begin{equation}
\bphi_j(\bx) = \cP^T \bphi_{j+1}(\bx), \qquad j \ge 0.
\label{eq:LS45}
\end{equation}
The Galerkin approximation of this equation,  in the same space used before, has the snapshots
\begin{equation}
\bphi_j^{\Gal,-}(\bx) = \bPhi(\bx)  \bg_j^{-}, \qquad j \ge 0,
\label{eq:LS45.1}
\end{equation}
where the  superscript $``-"$ reminds us that the stepping is backward in time and the coefficient matrices $\bg_j^-$ satisfy, similar to equation \eqref{eq:LS20}, 
\begin{equation}
\int_{\Omega} d \bx \, \bV^T(\bx) \left[ \bPhi(\bx)  \bg_j^{-} - \cP^T \bPhi(\bx)  \bg_{j+1}^{-} \right] = 
\bR  \bg_j^{-} - ({\cPR})^T \bR \bg_{j+1}^{-} = {\bf 0}, \qquad j \ge 0.
\label{eq:LS45.2}
\end{equation}
Here we used the Gram-Schmidt equation \eqref{eq:LS16} and the identity 
\begin{equation}
\int_{\Omega} d \bx \, \bV^T(\bx) \cP^T \bV(\bx) = \left[\int_{\Omega} d \bx \, \bV^T(\bx) \cP \bV(\bx) \right]^T 
\stackrel{\eqref{eq:LS33}} = ({\cPR})^T.
\end{equation}
The ROM snapshots for the backward iteration are  
\begin{equation}
\bphi_j^{\RM,-} = \bR  \bg_j^{-}, \qquad j \ge 0,
\label{eq:LS45.3}
\end{equation}
and they evolve according to the transpose of the ROM propagator matrix
\begin{equation}
 \bphi_j^{\RM,-} = ({\cPR})^T\bphi_{j+1}^{\RM,-} \qquad j \ge 0.
\label{eq:LS45.4}
\end{equation}
This is the (ROM) algebraic form of equation \eqref{eq:LS45}. 

Note that the ROM snapshots for the forward and backward iterations 
coincide at the first time instants because the same reasoning that lead to \eqref{eq:LS24} applies to the Galerkin approximation 
\eqref{eq:LS45.1}. We get 
\begin{equation}
\bg_j^{-} = \bg_j = \bi_j \quad \leadsto \quad \bphi_j^\RM =\bphi_j^{\RM,-} \qquad  \qquad j = 0, \ldots, n_t-1,
\label{eq:LS45.5}
\end{equation}
so the ROM propagator behaves as a unitary matrix for the first $n_t$ snapshots. More explicitly, if we substitute equation \eqref{eq:LS30} into \eqref{eq:LS45.4} and use equation \eqref{eq:LS45.5}, we get 
\begin{equation*}
\left[ \bI_{n_t n_\cE} - {\cPR}^T \cPR \right]\hspace{-0.07in} \begin{pmatrix} \bR_{0,0} & \bR_{0,1} & \bR_{0,2} & \ldots & \bR_{0,n_t-2} \\
{\bf 0} & \bR_{1,1} & \bR_{1,2} &  \ldots & \bR_{1,n_t-2} \\ 
&& \vdots && \\
{ \bf 0}& {\bf 0} & \ldots &{\bf 0} & \bR_{n_t-2,n_t-2} \\
 {\bf 0} & {\bf 0} & \ldots &{\bf 0} & {\bf 0} \end{pmatrix} = {\bf 0} \in \RR^{n_t n_{\cE} \times (n_t-1) n_{\cE}}.
\end{equation*}
Because the diagonal blocks of $\bR$ are invertible,  we deduce from this equation  that the first $n_{t}-1$ columns of $\big({\cPR}\big)^T \cPR$ equal those of the identity matrix.

\section{Acoustic waves}
\label{sect:Acoustics}

\begin{figure}[h]

\centering
\includegraphics[width=0.5\textwidth]{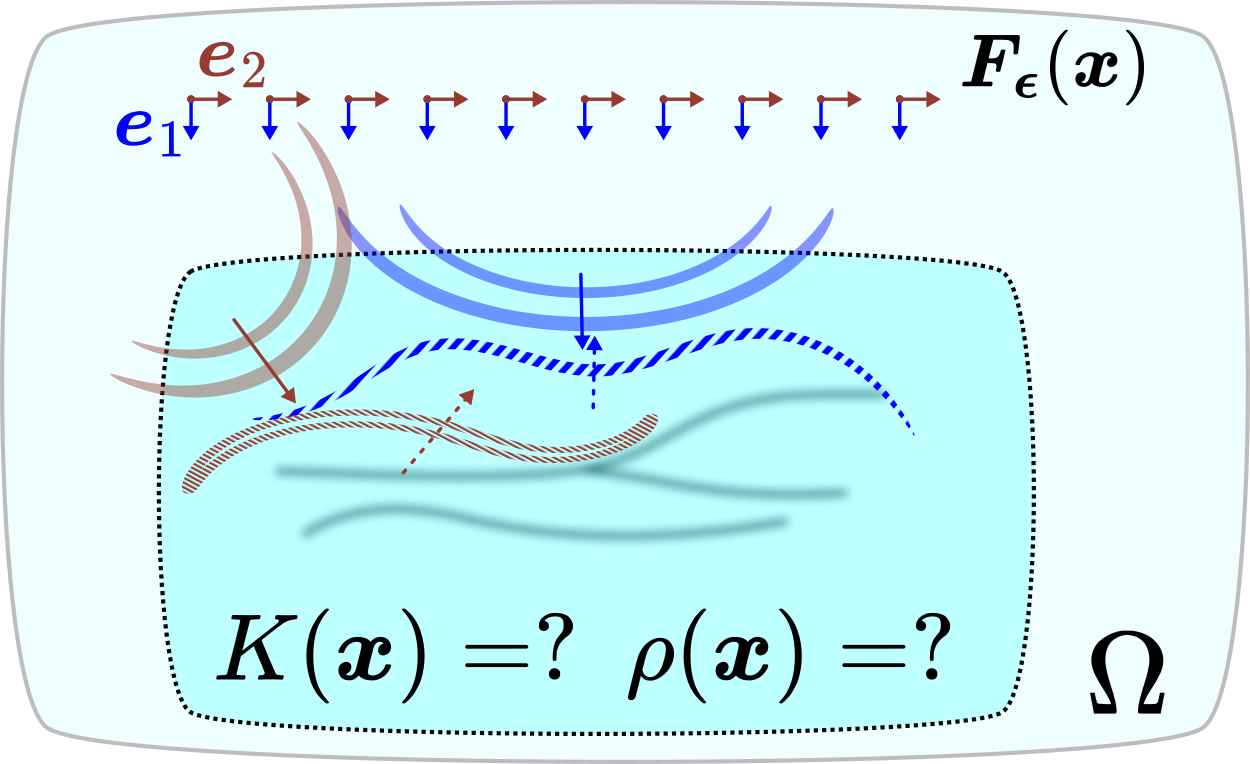}
\caption{Illustration of the acoustic setup for inverse scattering with polarized exitations ${\bm F}_{\bm \epsilon}( \bm x)$. The inversion domain is shown in 
darker blue and is surrounded by a homogeneous, reference medium. The sensors are placed above this domain and generate waves with $d$ different 
polarizations (here $d=2$).}
\label{fig:Drawing}
\end{figure}

In this section we specialize the ROM to acoustic waves, governed by the first order hyperbolic system
\begin{align}
\partial_t  \bu_{\be}(t,\bx) + \rho^{-1}(\bx) \nabla p_{\be}(t,\bx) &= \frac{s(t)}{\sqrt{\zeta_o}} \bF_{\be}(\bx)  , \label{eq:Ac1} \\
\partial_t p_{\be} (t,\bx) +  K(\bx) \nabla \cdot \bu_{\be} (t,\bx)  &= 0,
\label{eq:Ac2}
\end{align}
satisfied at $\bx \in \Omega$ and $t \in \RR$, by the acoustic pressure $p_{\be}$ and velocity $\bu_{\be}$, indexed by $\be = (\eps_1,\eps_2) \in \mathbb{N}^2$ (see Figure \ref{fig:Drawing}). Here the medium is modeled by the 
variable and unknown density $\rho$ and bulk modulus $K$, which determine the wave speed  
\begin{equation}
c(\bx) = \sqrt{\frac{K(\bx)}{\rho(\bx)}},
\label{eq:Ac3}
\end{equation}
and the wave impedance 
\begin{equation}
\zeta(\bx) = \sqrt{K(\bx)\rho(\bx)}.
\end{equation}
The reference values of the  coefficients are  the constants $\rho_o, K_o, c_o$ and $\zeta_o$, respectively. These model the 
homogeneous medium near the sensors in the array.

\subsection{Derivation of the first order system \eqref{eq:LS1}} 
\bc{The theory in the previous section was done for arbitrary forcing. For its application to acoustics, we assume a separable 
force function, which is typical in phased array data acquisitions. Thus, the right hand side in equation \eqref{eq:Ac1} consists of the  pulsed signal 
$s$, supported at $t \in (-t_s,t_s)$, and }
\begin{equation}
\bF_{\be}(\bx) \approx \delta(\bx-\bx_{\eps_1}){\itbf e}_{\eps_2},
\label{eq:Ac4}
\end{equation}
where $\bx_{\eps_1}$ is the location of the $\eps_1^{\rm th}$ sensor  and ${\itbf e}_1, \ldots, {\itbf e}_d$ is the canonical orthonormal basis in $\RR^d$. The approximation means 
that $\bF_{\be}$ is an $\bL^2(\Omega)$ function normalized by 
\begin{equation}
\int_{\Omega} d \bx \, \bF_{\be}(\bx) = {\itbf e}_{\eps_2},
\label{eq:Ac5}
\end{equation}
and supported in a ball centered at  $\bx_{\eps_1}$, with small radius with respect to the wavelength.
The scaling of the force by ${\zeta_o}^{-1/2}$   in \eqref{eq:Ac1} is to simplify the formulas below.
We model the boundary of the domain as sound soft, 
\begin{equation}
p_{\be}(t,\bx) = 0, \qquad t \in \RR, ~~ \bx \in \partial \Omega,
\label{eq:Ac6}
\end{equation}
and suppose that 
the medium is quiescent prior to the excitation, 
\begin{equation}
p_{\be}(t,\bx) = 0, \qquad \bu_{\be}(t,\bx) = {\bf 0}, \qquad t \ll 0, ~~ \bx \in \Omega.
\label{eq:Ac7}
\end{equation}

The hyperbolic system \eqref{eq:Ac1}-\eqref{eq:Ac2} can be put in the desired form \eqref{eq:LS1}, by defining the 
following $m = d+1$ dimensional wave field
\begin{equation}
\bpsi_{\be}(t,\bx) = \sqrt{c_o} \begin{pmatrix} 
\sqrt{\rho(\bx)} \bu_\eps(t,\bx) \\ 
\frac{1}{\sqrt{K(\bx)}} p_\eps(t,\bx) 
\end{pmatrix}.
\label{eq:Ac8}
\end{equation}
We obtain equation \eqref{eq:LS1} with the skew-adjoint operator
\begin{equation}
\cL = \begin{pmatrix} 0 & \frac{1}{\sqrt{\rho(\bx)}} \mbox{grad} \left[ c(\bx) \sqrt{\rho(\bx)} \cdot \right] \\
c(\bx) \sqrt{\rho(\bx)} \mbox{div} \left[ \frac{1}{\sqrt{\rho(\bx)}} \cdot \right] & 0 \end{pmatrix} ,
\label{eq:Ac9}
\end{equation}
and the force
\begin{equation}
\bc{ 
\bff_{\be}(t,\bx) = s(t) \begin{pmatrix}\bF_{\be}(\bx) \\ 0 \end{pmatrix}.
}
\label{eq:Ac10}
\end{equation}

\subsection{The data} The variable coefficients in $\cL$, written in terms of the wave speed $c$ and density $\rho$, are the unknowns, to be determined 
from the measured array response matrix $\bA$. 
The entries  \eqref{eq:LS3} of this matrix are 
\begin{align}
\nonumber
\bc{\bA_{\be',\be}(t)} &\bc{\stackrel{\eqref{eq:Ac10}}{=} \int_{-t_s}^{t_s} dt' \int_{\Omega} d \bx \, s(-t') \begin{pmatrix}\bF_{\be'}(\bx) \\ 0 \end{pmatrix}^T \bpsi_{\be}(t-t',\bx) } \\
\nonumber
& 
\bc{ \stackrel{\eqref{eq:Ac8}}{=}  \int_{-t_s}^{t_s} dt' \int_\Omega d\bx \sqrt{c_o} \sqrt{\rho(\bx)}  s(-t') \bF_{\be'}^T(\bx) \bu_{\be}(t-t',\bx ),} \\
&
\bc{ \approx  \int_{-t_s}^{t_s} dt' \sqrt{\zeta_o} {\itbf e}_{\eps'_2}^T s(-t')\bu_{\be}(t-t',\bx_{\eps'_1}),}
\label{eq:Ac11}
\end{align}
where the approximation is due to equation \eqref{eq:Ac4} and the assumption $\rho(\bx_{\eps'_1}) = \rho_o$. Thus, the array measures the velocity at the sensor locations. These measurements  can be related to measurements of the pressure using the divergence theorem,
\begin{align}
\nonumber&\bc{
\int_{-t_s}^{t_s} dt' \int_{\Omega} d \bx \, s(-t') \mbox{div}[ \bF_{\be'}(\bx)] p_{\be}(t-t',\bx)}\\
 &\bc{= 
 \int_{-t_s}^{t_s} dt'  s(-t') 
\int_{\Omega} d \bx \, \mbox{div}[\bF_{\be'}(\bx)p_{\be}(t-t',\bx) ]  }
\nonumber\\
&\quad \bc{-  \int_{-t_s}^{t_s} dt'  \int_{\Omega} d \bx \left[\bF_{\be'}(\bx)\right]^T  s(-t') \nabla p_{\be}(t-t',\bx)}\nonumber \\
&\hspace{-0.06in}\bc{\stackrel{\eqref{eq:Ac1}}{=} -  \int_{-t_s}^{t_s} dt' \int_{\Omega} d \bx \left[ \bF_{\be'}(\bx) \right]^T    s(-t') \left[ \frac{s(t-t') \rho_o}{\sqrt{\zeta_o}} \bF_{\be}(\bx) - 
\rho(\bx) \partial_t \bu_{\be}(t-t',\bx) \right] }\nonumber \\
&\bc{= \frac{\rho_o}{\sqrt{\zeta_o}} \left[\frac{d}{dt} \bA_{\be',\be}(t) - \int_{-t_s}^{t_s} dt'   s(-t') s(t-t') \int_{\Omega} d \bx \left[ \bF_{\be'}(\bx) \right]^T \bF_{\be}(\bx)\right].}
\label{eq:Ac12}
 \end{align}
 Here we used the boundary condition \eqref{eq:Ac6} and that $\rho = \rho_o$ at the array. Note that the last term in the right hand side
of this equation vanishes when $\eps_1' \ne \eps_1$, because the two forces have disjoint support. Note also that given our model \eqref{eq:LS5} 
of the force, the derivative of the array response matrix gives approximately the components of the gradient of pressure at the array
\begin{align}
\nonumber
&
\bc{\int_{-t_s}^{t_s} dt'
\int_\Omega d \bx \,  s(-t') \mbox{div}[ \bF_{\be'}(\bx)] p_{\be}(t-t',\bx)}\\
\nonumber
& \bc{\stackrel{\eqref{eq:Ac4} }{\approx} \int_\Omega d \bx \,  \partial_{x_{\eps'_2}} \delta(\bx-\bx_{\eps'_1}) \int_{-t_s}^{t_s} dt' s(-t') p_{\be}(t-t',\bx)}\\
&
\bc{= -   \int_{-t_s}^{t_s} dt' s(-t') \partial_{x_{\eps'_2}} p_{\be}(t-t',\bx_{\eps'_1})  .}
\label{eq:Ac13}
\end{align}
Suppose that the array is planar and that we measure the components of the velocity in the normal direction to the array. The calculation above says 
that this is basically the same as measuring the normal derivative of the pressure at the array. This is not enough information for solving the inverse problem, so we need to measure the velocity in the plane of the array, as well.

We end this section with the proof of reciprocity (symmetry of $\bA$), mentioned in Remark \ref{rem:rec}: Equation \eqref{eq:Ac12} shows that $\bA$ is symmetric, 
if the pressure field satisfies
\begin{equation}
\int_{\Omega} d \bx \, \mbox{div}[\bF_{\be'}(\bx)] p_{\be}(t,\bx) = \int_{\Omega} d \bx \, \mbox{div}[\bF_{\be}(\bx)] p_{\be'}(t,\bx),  \qquad 
\forall ~\be, \be' \in \cE.
\label{eq:Rec1}
\end{equation}
Equations \eqref{eq:Ac1}-\eqref{eq:Ac2} give that $p_\be$ satisfies the second order wave equation 
\begin{equation}
\partial_t^2 p_{\be}(t,\bx) - K(\bx) \mbox{div} \left[ \rho^{-1}(\bx) \nabla p_\be(t,\bx)\right] = - c_o \frac{d}{dt} s(t) \mbox{div}[\bF_{\be}(\bx)],
\end{equation}
with homogeneous Dirichlet boundary condition at $\partial \Omega$. If we let $G$ be the Green's function for the wave operator in 
this equation, we get that  
\begin{equation}
p_{\be}(t,\bx) = c_o \frac{d}{dt} s(t) \star_t \int_{\Omega}d \by \,  G(t,\bx,\by)  \mbox{div}[\bF_{\be}(\by)],
\label{eq:Rec1p}
\end{equation}
where $\star_t$ denotes convolution in time. It is known that $G$ satisfies the reciprocity relation
$
G(t,\bx,\by) = G(t,\by,\bx)
$
(see for example  \cite{fokkema2013seismic}). Therefore, equation \eqref{eq:Rec1p} implies that \eqref{eq:Rec1} holds.

\section{Multiparametric inversion}
\label{sect:Inversion}
If we knew the wave field inside the domain $\Omega$, it would be easy to estimate $c$ and $\rho$. This is precisely the point of hybrid (multi-physics) imaging modalities \cite{bal2013hybrid} that use one kind of wave to probe the medium and then measure the 
propagation of another wave. Hybrid methods require complex equipment and controlled measurement 
settings that cannot be implemented in general environments. Here we use the  ROM to map the measurements at the array to a good approximation of the 
vectorial wave field inside $\Omega$. 

The basic idea for obtaining such a map comes from the  Gram-Schmidt equation \eqref{eq:LS16}, that  factorizes the 
unknown snapshots contained in $\bPhi$ in the data driven ROM snaphsots, the block columns of $\bR$, and the uncomputable 
orthonormal basis  in $\bV$. We know from Theorem \ref{thm.1} that the ROM snapshots satisfy the exact 
data fitting relation \eqref{eq:LS34}, just like the true snapshots. This means that all the information registered at the 
array  is contained in the ROM snapshots. The uncomputable basis $\bV$ plays no role in the data fitting. Indeed, equations  \eqref{eq:LS16} 
\eqref{eq:LS33} and \eqref{eq:LS40} give 
\begin{align}
\bD_j &= \int_\Omega d \bx \, \bphi_0^T(\bx) \bphi_j(\bx) \nonumber = \int_\Omega d \bx \, \left[\bV(\bx) \bphi_0^{\RM}\right]^T \bV(\bx) \bphi_j^{\RM}  \nonumber \\
&= 
\big(\bphi_0^{\RM}\big)^T \underbrace{\int_\Omega d \bx \, \bV^T(\bx) \bV(\bx)}_{\bI_{n_t n_\cE}}  \bphi_j^{\RM} = \big(\bphi_0^{{\RM}}\big)^T  \bphi_j^{\RM}, \qquad j = 0, \ldots, n_t-1, 
\label{eq:N1}
\end{align}
and the same result  holds when replacing $\bV$ with  any orthonormal basis. That is to say, if we define 
\begin{equation}
\tilde \bphi_j(\bx;\tilde c,\tilde \rho) = \bV(\bx;\tilde c,\tilde \rho) \bphi_j^{\RM}, \qquad j = 0, \ldots, n_t-1, 
\label{eq:N2}
\end{equation}
with $\bV(\cdot;\tilde c, \tilde \rho)$ the orthonormal basis computed via the Gram-Schmidt procedure in the medium with search wave speed $\tilde c$ and density $\tilde \rho$, 
we get, similar to \eqref{eq:N1}, the exact data fit
\begin{align}
\int_\Omega d \bx \, \tilde \bphi_0^T(\bx;\tilde c, \tilde \rho) \tilde\bphi_j(\bx;\tilde c,\tilde \rho) 
&= 
\big(\bphi_0^{{\RM}}\big)^T \underbrace{\int_\Omega d \bx \, \tilde \bV^T(\bx;\tilde c, \tilde \rho) \bV(\bx;\tilde c, \tilde \rho)}_{\bI_{n_t n_\cE}} \bphi_j^{\RM}  \nonumber \\
&=\big( \bphi_0^{{\RM}}\big)^T  \bphi_j^{\RM} = \bD_j, \qquad j = 0, \ldots, n_t-1.
\label{eq:N3}
\end{align}

The basis $\bV$ maps the information registered at the array and contained in the ROM snapshots, from the algebraic (ROM) space to 
the physical space, at  $\bx \in \Omega$. Its computation requires knowing the medium. Thus, we define our approximation of the 
internal wave by \eqref{eq:N2}, where the uncomputable $\bV$ is replaced by the basis in our best guess of the medium, modeled 
by $\tilde c$ and $\tilde \rho$. While this internal wave is consistent with the data,  as shown above, it is not a solution of the hyperbolic system of 
equations. Assuming that the inverse problem is uniquely solvable \cite{bhattacharyya2022recovery}, this  would be the case only when $\tilde c$ and $\tilde \rho$ are equal to the true wave speed and density. 

Therefore, our inversion approach seeks to minimize the discrepancy ${\cal O}(\tilde c, \tilde \rho)$ between the estimated internal waves \eqref{eq:N2} and the solution of the wave equation, whose first $n_t$ snapshots are gathered in $\bV(\cdot;\tilde c, \tilde \rho) \bR(\tilde c, \tilde \rho)$, similar to equation \eqref{eq:LS16}. 
The objective function is 
\begin{align}
 {\cal O}(\tilde c, \tilde \rho) &= \sum_{j=0}^{n_t-1} \int_{\Omega} d \bx \, \|\tilde \bphi_j(\bx;\tilde c,\tilde \rho) - \bphi_j(\bx;\tilde c,\tilde \rho) \|_F^2 \nonumber \\
 &= \int_{\Omega} d \bx \, \|\bV(\bx;\tilde c, \tilde \rho) \bR - \bV(\bx;\tilde c, \tilde \rho) \bR(\tilde c, \tilde \rho)\|_F^2  
\nonumber \\
 &= \|\bR -  \bR(\tilde c, \tilde \rho)\|_F^2,
 \label{eq:N4}
 \end{align}
 where $\| \cdot \|_F$ is the Frobenius norm,  $\bR$ without an argument is the data driven matrix containing the ROM snapshots  and $\bR(\tilde c, \tilde \rho)$ is the block Cholesky square root of the mass matrix computed in the guess medium 
 with wave speed $\tilde c$ and density $\tilde \rho$. The contribution of $\bV(\cdot;\tilde c, \tilde \rho)$ cancels out in \eqref{eq:N4} by the definition of the Frobenius norm  and orthonormality. Thus, we are left with minimizing the misfit of the block Cholesky square root of the mass matrix.

Instead of working with  \eqref{eq:N4}, we minimize the objective function 
\begin{align}
 {\cal O}^{\RM} (\tilde c, \tilde \rho) &= \| \bR(\tilde c, \tilde \rho) \bR^{-1}-\bI_{n_t n_\cE} \|_F^2.
 \label{eq:N5}
 \end{align}
 This works better because by taking the inverse of $\bR$ we balance  the contribution of the strong and weak echoes 
 registered in the data $\{\bD_j\}_{j=0}^{n_t-1}$.

\subsection{The inversion algorithm}
\label{sect:alg} 
To formulate the algorithm, we need to specify the search space where $\tilde c$ and $\tilde \rho$ lie. \bc{We  
parameterize $\tilde c$ and $\tilde \rho$ as}
\begin{align}
\tilde c(\bx; \bet) & = c_o + \sum_{l=1}^{N} \eta_l \beta_l(\bx), 
\qquad 
\tilde \rho (\bx; \bet)  = \rho_o + \sum_{l=1}^{N} \eta_{l+N} \beta_l(\bx),
\label{eqn:rparam}
\end{align}
where we recall that $c_o$ and $\rho_o$ are the known reference coefficients. \bc{
Note that  in \eqref{eqn:rparam} we use the same basis functions $\{\beta_j\}_{j=1}^N$ for $\tilde c$ and $\tilde \rho$, because we expect that the variations of the
wave speed  and density are related, since they model the same heterogeneous medium. However, the inversion algorithm is not tied to using the same basis and the same number of parameters. One could use independent parameterizations for the two unknowns}. 

The objective function \eqref{eq:N5}, for the parameterization 
\eqref{eqn:rparam},
depends on the vector 
\[\bet = (\eta_1,\ldots,\eta_N,\eta_{N+1},\ldots,\eta_{2N})^T ,\] so we simplify the notation as 
\begin{align}
{\cal O}^{\RM} (\bet) = {\cal O}^{\RM} (\tilde c(\cdot;\bet), \tilde \rho(\cdot;\bet)), \qquad 
\bR(\bet) = \bR(\tilde c(\cdot;\bet), \tilde \rho(\cdot;\bet)).
\label{eq:N6}
\end{align}

\bc{The causality  of the ROM, ensured by the  Gram-Schmidt orthogonalization (\ref{eq:LS16}) and reflected in the causal (block upper triangular) structure of $\bR$, can be used to  estimate the medium in a  layer stripping fashion: 
Denote by $[\bR]_k$ the upper left sub-matrix of $\bR$, containing the first 
$k$ blocks. This is defined by the early time data $\bD_j$, $j=0,\ldots,k-1$, for $k < n_t$, which contain information about the medium 
sampled by the waves up to time $t_{k-1}$. Consequently, we can 
feed the data gradually, for increasing values of $k$, and estimate the medium in layers,  by minimizing the sub-objective functions
\begin{equation}
\mathcal{O}^\RM_k(\bet) = \left\| [\bR(\bet)]_k ([\bR]_k)^{-1} - \bI_{k n_{\cE}}  \right\|_F^2.
\label{eqn:ok}
\end{equation}
}

As is the case with any inversion algorithm, the objective function must be complemented by a regularization penalty term. Any such regularization would work. We use the simplest,  Tikhonov type penalty,
\begin{equation}
\mathcal{O}^\REG(\bet) = \| \bet \|_2^2.
\label{eqn:oreg}
\end{equation}
The minimization of the objective function \eqref{eqn:ok} complemented by the  regularization penalty can be carried out with any gradient based method. We choose the Gauss-Newton method.

\begin{algorithm}\textbf{\emph{(ROM based inversion)}}
\label{alg:rominv}

\vspace{0.04in} \noindent  \textbf{Input:} The data driven $\bR$; the number $\ell$ of layers for layer stripping 
and the number $q$ of iterations per layer; integer layer indices $\{k_l\}_{l=1}^\ell$, 
such that \[1 \le k_1 \leq k_2 \leq  \cdots \leq  k_{\ell} = n_t.\]

\vspace{0.04in} \noindent  Starting with the initial vector $\bet^{(0)}= {\bf 0}$, proceed:\\ 
For $l = 1,2,\ldots, \ell$, and $j = 1,\ldots, q$, set the update index \[i = (l-1)q+j.\]
Compute $\bet^{(i)}$ as a Gauss-Newton update for minimizing the functional
\begin{equation}
\mathcal{R}_i^\RM(\bet) = \mathcal{O}^\RM_{k_l}(\bet) + \mu_i \mathcal{O}^\REG(\bet),
\label{eqn:objreg}
\end{equation}
linearized about $\bet^{(i-1)}$. The penalty weight $\mu_i > 0$ is chosen adaptively at each iteration,
based on the decay of the singular values
of the Jacobian of $\mathcal{R}_i^\RM(\bet^{(i-1)})$. See \cite{borcea2022waveform} for the detailed 
description of the penalty weight computation.

\vspace{0.04in} \noindent  \textbf{Output:} 
The velocity and density estimates computed from equations \eqref{eqn:rparam}, with the coefficients in 
$\bet^{(\ell q)}$.
\end{algorithm}

\subsection{Regularization of the ROM computation}
\label{sec:romreg}
So far we  supposed that Assumption \ref{as:1} holds and thus,  the mass matrix $\bbM$ is symmetric and positive definite. 
In reality this may not be true, especially for noisy data, so $\bbM$ must be regularized before we can compute its block Cholesky square root. 

It is easy to symmetrize the data driven mass matrix $\bbM$, by taking 
$(\bbM + \bbM^T)/2$. 
\bc{This symmetrization is, in fact, the maximal likelihood estimator of a symmetric matrix given a noisy observation with additive Gaussian independent noise.}
We assume henceforth that $\bbM$ is symmetric so it remains to map it to a positive definite matrix. This must be 
done carefully, in order to preserve the causality of the ROM, which is essential for a successful inversion. Here is one way to 
do it: Let 
\begin{equation}
\bbM = {\itbf Z} \boldsymbol{\Lambda}{\itbf Z}^T,
\label{eqn:eigenmass}
\end{equation}
be the eigendecomposition of $\bbM$, where ${\itbf Z}$ is the orthogonal matrix of 
eigenvectors and $\boldsymbol{\Lambda} = {\rm diag} \big(\la_1, \ldots, \la_{n_t n_{\cE}}\big)$ 
is the diagonal matrix of eigenvalues, sorted in descending order. We wish to filter out the 
smallest eigenvalues and the associated eigenvectors. Since the mass matrix has a block structure with $n_{\cE} \times n_{\cE}$ blocks, 
we choose a number $r$ such that $1 \le r < n_t$, and  keep the $r n_{\cE}$ largest eigenvalues of $\bbM$. 
The corresponding leading orthonormal eigenvectors are  the columns of
\begin{equation}
{\itbf Z}^{(r)}
 = [Z_{jk}]_{1\leq j \leq n_t n_{\cE}, 1 \leq k \leq r n_{\cE}}
\in \RR^{n_t n_{\cE} \times r n_{\cE}}
\label{eqn:zr} 
\end{equation}
and the projected mass matrix is
\begin{equation}
\boldsymbol{\Lambda}^{(r)} = ( {\itbf Z}^{(r)})^T \bbM {\itbf Z}^{(r)}
 = \mbox{diag} \big(\la_1, \ldots, \la_{r n_{\cE}}\big).
\label{eqn:mproj}
\end{equation}
The stiffness matrix ${\bbS}$ can also be projected
\begin{equation}
{\bbS}^{(r)} = ( {\itbf Z}^{(r)})^T {\bbS} {\itbf Z}^{(r)} 
\in \RR^{r n_{\cE} \times r n_{\cE}},
\label{eqn:sntilde}
\end{equation}
and  the analogue of the ROM propagator defined in \eqref{eq:LS31} is
\begin{equation}
\mathbb{P}^{(r)} = (\boldsymbol{\Lambda}^{(r)})^{-1/2} 
{\bbS}^{(r)} (\boldsymbol{\Lambda}^{(r)})^{-1/2} \in \RR^{r n_{\cE} \times r n_{\cE}}.
\label{eqn:pproj}
\end{equation}

For a proper choice of $r$ we have a well conditioned projected mass matrix \eqref{eqn:mproj},  at the cost of losing the block Toeplitz structure in 
\eqref{eq:LS36}--\eqref{eq:LS37} and the block upper Hessenberg structure of the ROM propagator. These structures reflect the 
causality of the ROM. To regain them, we use an additional orthogonal transformation, given by  the  block Arnoldi Algorithm \cite{sadkane1993block} applied to  the matrix $\mathbb{P}^{(r)}$, with starting block
\footnote{The published version of the paper 
\href{https://doi.org/10.1137/24M1699784}{DOI:10.1137/24M1699784} 
contains a typo in \eqref{eqn:startblock} that is corrected here.} 
\begin{equation}
(\boldsymbol{\Lambda}^{(r)})^{1/2} ({\itbf Z}^{(r)})^T \bi^{(r)}_0 
\in \RR^{r n_{\cE} \times n_{\cE}},
\label{eqn:startblock}
\end{equation}
where $\bi^{(r)}_0 \in \RR^{r n_{\cE} \times n_{\cE}}$ contains the first $n_{\cE}$ columns of the 
$r n_{\cE} \times r n_{\cE}$ identity matrix. We give the algorithm in  Appendix \ref{ap:B}. It produces an orthogonal matrix 
${\itbf Q}^{(r)}$ such that 
$({\itbf Q}^{(r)})^T \mathbb{P}^{(r)} {\itbf Q}^{(r)}$ is
block upper Hessenberg. Moreover, ${\itbf Q}^{(r)}$
enables the computation of the block Cholesky square root
\begin{equation}
({\itbf Q}^{(r)})^T \boldsymbol{\Lambda}^{(r)} {\itbf Q}^{(r)}
= (\boldsymbol{\Pi}^{(r)})^T \bbM \boldsymbol{\Pi}^{(r)} 
=  (\bR^{(r)})^T \bR^{(r)}, \qquad \boldsymbol{\Pi}^{(r)} = {\itbf Z}^{(r)} {\itbf Q}^{(r)}.
\label{eqn:rnoise}
\end{equation}

 In inversion, we also need to compute the block Cholesky square
root of the mass matrix $\bbM(\bet)$ for the current search velocity and density parameterized as in \eqref{eqn:rparam}. 
We make this calculation compatible with the regularization process 
\eqref{eqn:eigenmass}--\eqref{eqn:rnoise}, by using the same $\boldsymbol{\Pi}^{(r)}$ defined in \eqref{eqn:rnoise} to get 
\begin{equation}
(\boldsymbol{\Pi}^{(r)})^T  \bbM(\bet)  \boldsymbol{\Pi}^{(r)} = 
\bR^{(r)}(\bet)^T \bR^{(r)}(\bet).
\label{eqn:rvrnoise}
\end{equation}
The inversion Algorithm~\ref{alg:rominv} works exactly the same, except that $\bR$ and $\bR(\bet)$ are 
replaced by $\bR^{(r)}$ and $\bR^{(r)}(\bet)$ defined above. 


\subsection{Numerical results}
\label{sect:Numerics}
In this section we compare  the performance of Algorithm~\ref{alg:rominv} to conventional
multiparameter full waveform inversion (MFWI). To make the comparison as close as possible, our implementation of 
MFWI uses the same data matrices as Algorithm~\ref{alg:rominv}, but instead of \eqref{eqn:objreg}, it 
minimizes the objective function 
\begin{equation}
\mathcal{R}_i^\FWI(\bet) = \mathcal{O}^\FWI_{k_l}(\bet) + \mu_i \mathcal{O}^\REG(\bet), \qquad \mathcal{O}^\FWI_{k_l}(\bet) = \sum_{j=0}^{k_l} 
\| \text{Triu}\left( \bD_j - \bD_j(\bet) \right) \|_F^2,
\label{eqn:objmfwireg}
\end{equation}
where $\bD_j(\bet)$ are the data matrices computed in the medium with search speed $\tilde c$ and  density $\tilde \rho$, 
parameterized by $\bet$ as in equation \eqref{eqn:rparam}. Since the noiseless matrices $\bD_j$ are symmetric, we use the operation 
$\text{Triu}$ to zero out the entries below the main diagonal.

We refer to Appendix \ref{ap:A} for the details of the numerical simulations, the model of noise  added to the data, \bc{and the choice of the regularization parameter $r$}. The remainder 
of the section is dedicated to presenting numerical results. In all the simulations we start the optimization with  $\bet = {\bf 0}$. 
{The number of iterations for the MFWI and ROM algorithms is chosen to be the same, so that for both of them, the last updates are of the order of $10^{-2} \| c_0 \|_\infty$ and $10^{-2} \| \rho_0 \|_\infty$ of the sup norm $\| \cdot \|_{\infty}$ of the wave speed and density, respectively.}

\bc{From various numerical experiments, we observed that  the time step $\Delta t$, the aperture size and the number of sensors affect the quality of the inversion. 
The effect of $\Delta t$ and the sensor separation  on the approximation of the internal wave  was already investigated in the scalar acoustic case in \cite{borcea2022internal}, and it is similar here. Roughly, the results say that the optimal choice is for $\Delta t \lesssim 1/(2 \nu)$ and sensor separation of the order of half the wavelength. Here, we make the additional observation  that  the estimation of the density benefits from having a large aperture. The choices made in the experiments below follow these observations.}

\subsection{Disjoint inclusions}
The objective of the first numerical experiment is to assess the amount of cross-talk between the 
two acoustic parameters. The medium has two disjoint rectangular inclusions, one in $c$ and 
another in $\rho$ in a homogeneous background with $c_o = 1$km/s and $\rho_o = 1\mbox{g/cm}^3$. 
The inclusion contrasts are $1.5$ for both $c$ and $\rho$. The domain is  $\Omega = [0, 3\mbox{km}] \times [0, 2\mbox{km}]$. 
The array has  $n_s=35$ sensors at depth of $1\mbox{km}$, chosen to minimize the effect of the boundary on this experiment. The data are noiseless.  The optimization search space 
has dimension $2 N = 2(30 \times 20) = 1200$ with velocity and density parameterized as in 
\eqref{eqn:rparam}.

\begin{figure}[h]
\begin{center}
\begin{tabular}{cc}
$c\; (km/s)$ & $\rho\; (g/cm^3)$ \\
\includegraphics[width=0.47\textwidth]
{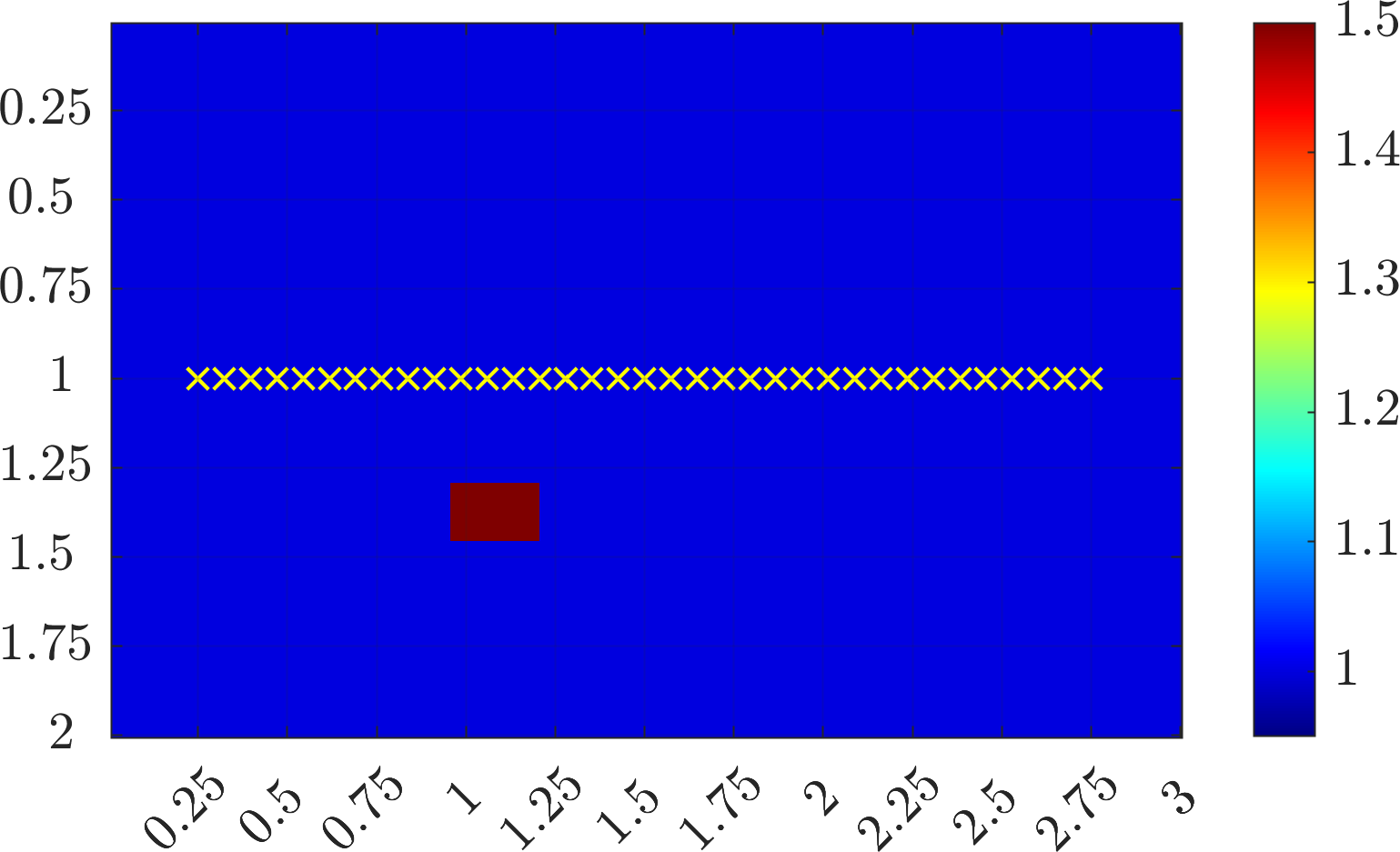} &
\includegraphics[width=0.47\textwidth]
{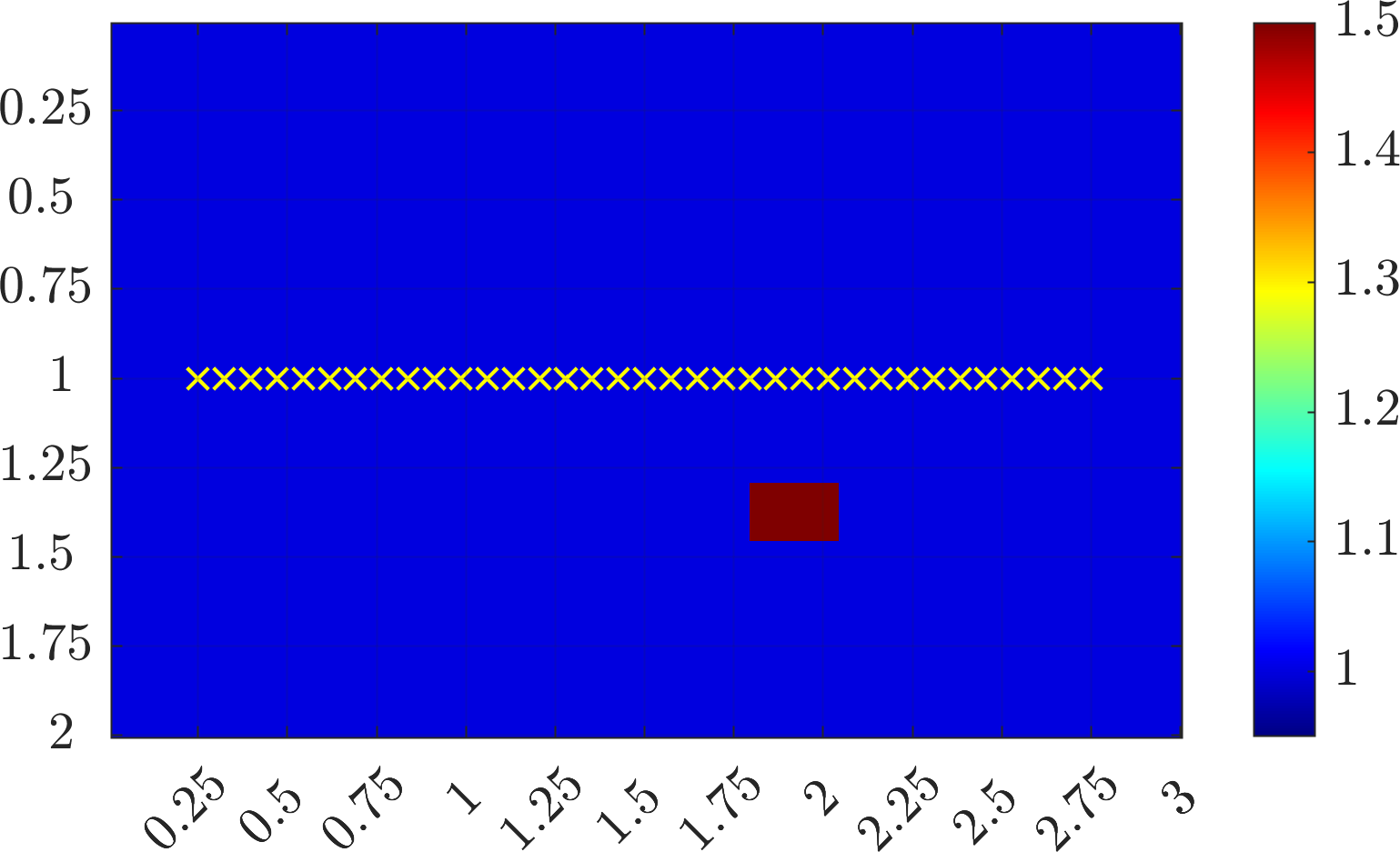} \\
\includegraphics[width=0.47\textwidth]
{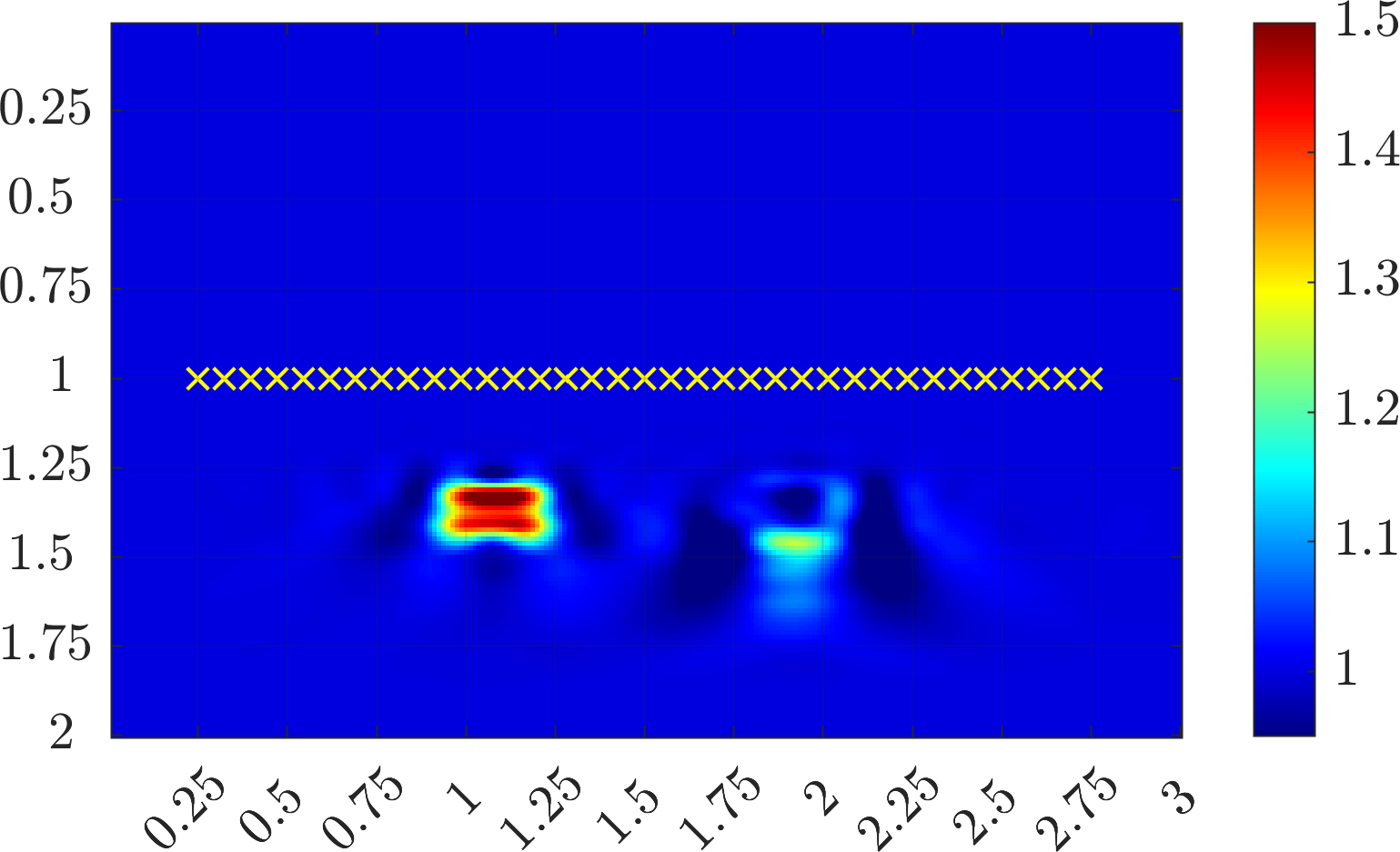} &
\includegraphics[width=0.47\textwidth]
{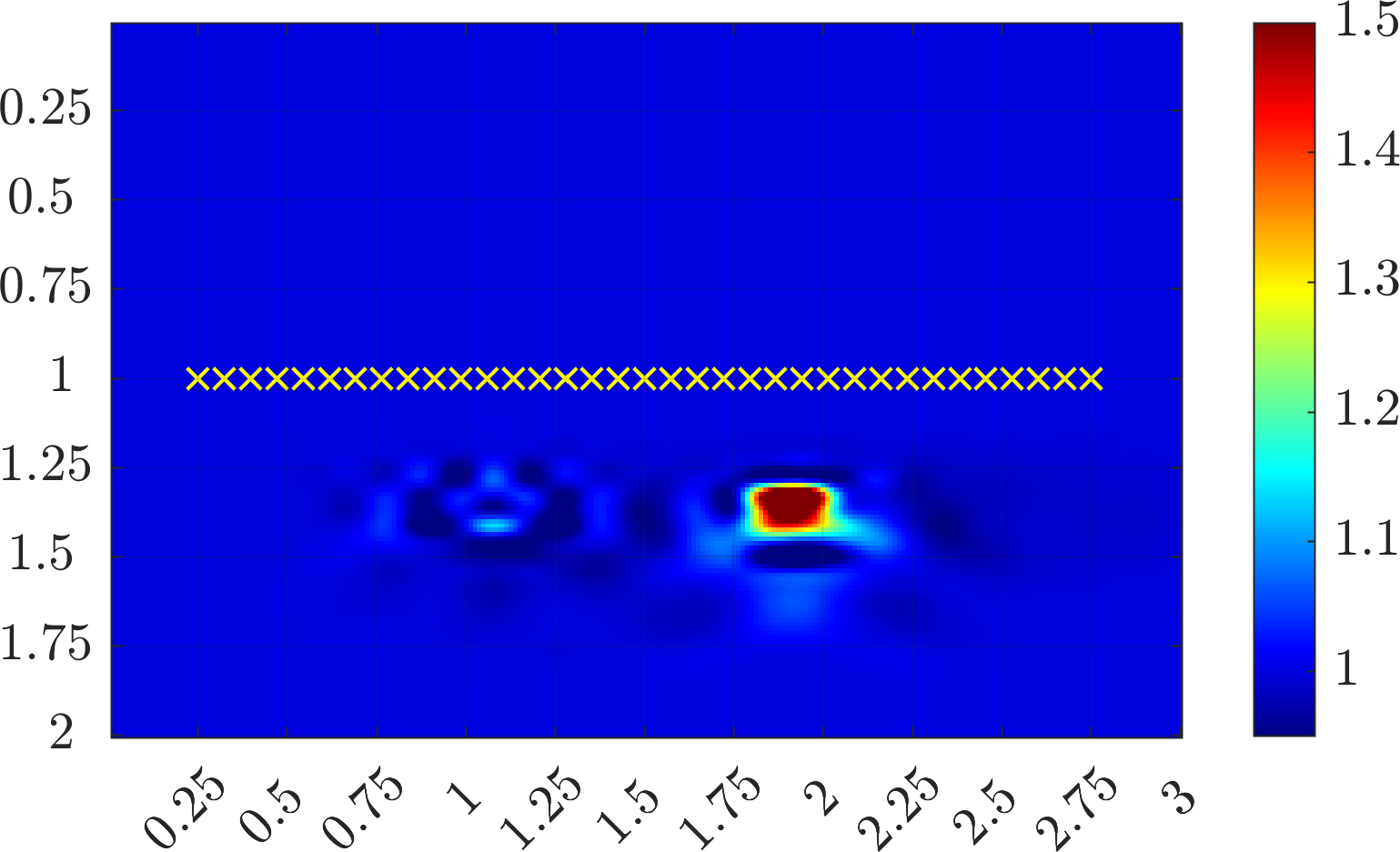} \\
\includegraphics[width=0.47\textwidth]
{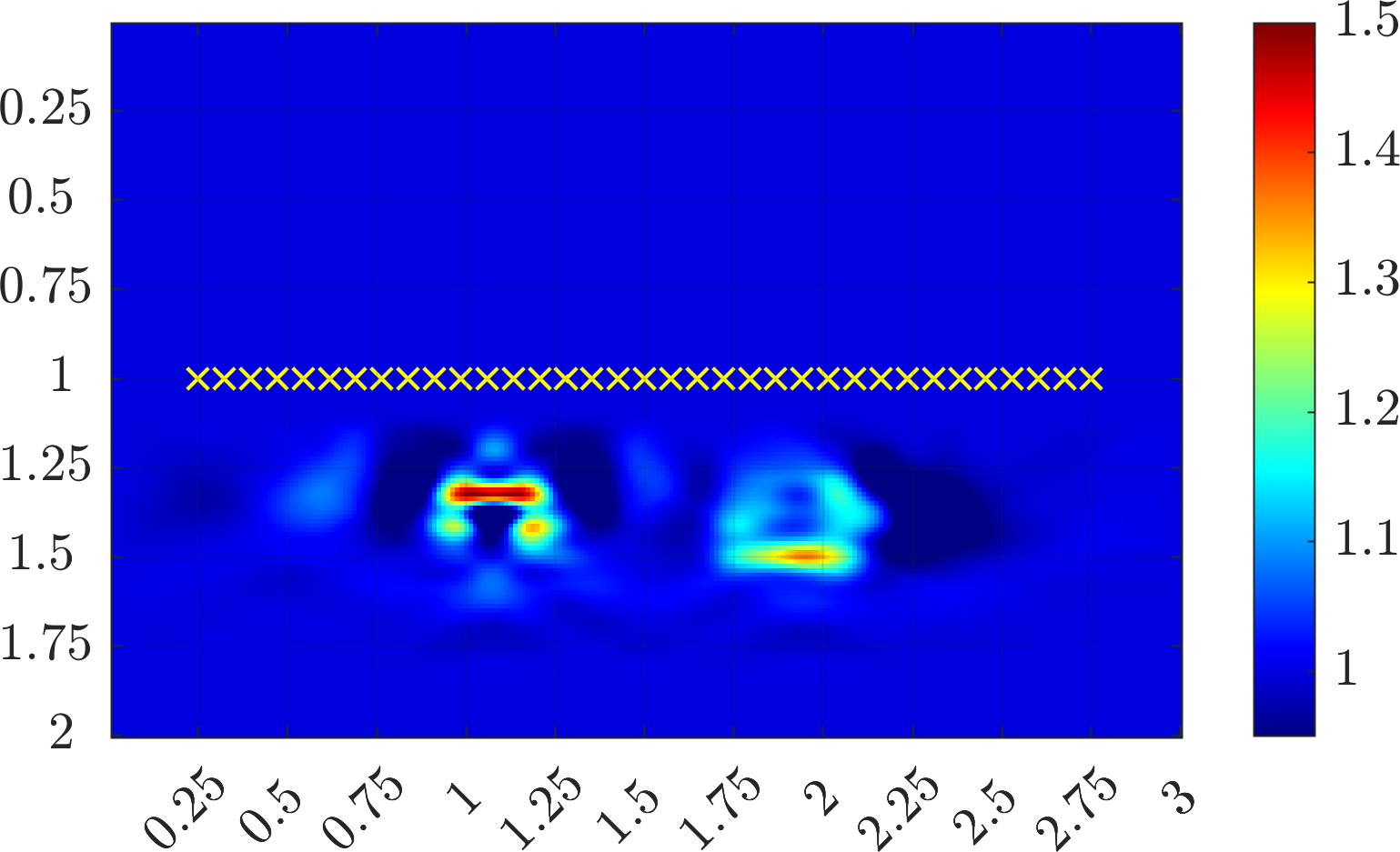} &
\includegraphics[width=0.47\textwidth]
{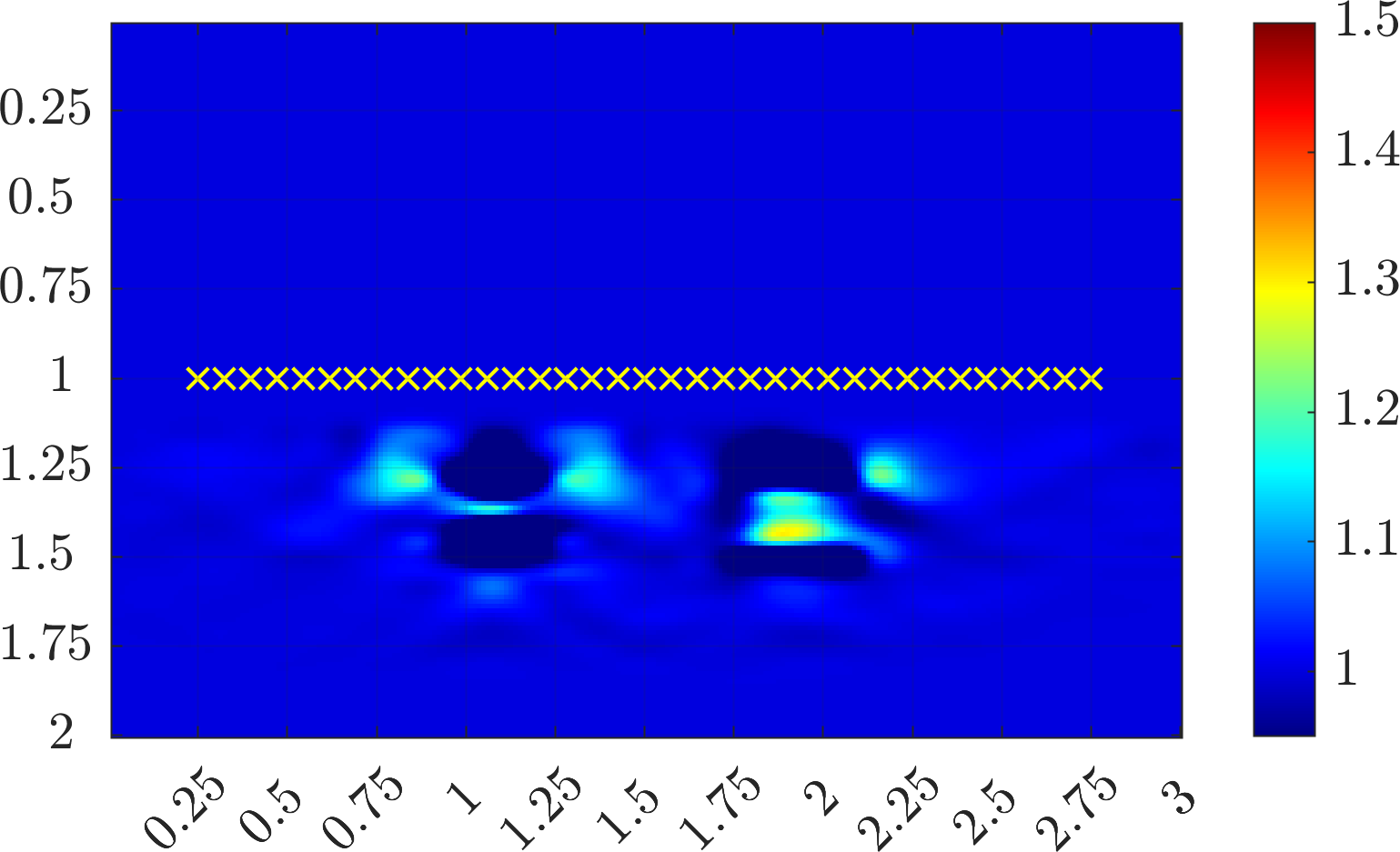}
\end{tabular}
\end{center}
\caption{Two disjoint inclusions. Top row: true velocity and density; 
middle row: ROM based estimates; bottom row: MFWI estimates.
Source locations are yellow crosses. Axes are in km. The colorbar shows the contrast i.e., 
ratios of $c$ and $\rho$ with the reference values $c_o$ and $\rho_o$ of the medium near the array.
}
\label{fig:2inc}
\end{figure}

We display in Figure~\ref{fig:2inc} the model and the estimates 
obtained with MFWI and our ROM algorithm. For both approaches we take $15$ iterations. 
Given the simplicity of the model,
we do not use layer stripping in this experiment i.e., $\ell = 1$, $k_1 = n_t$, and $q = 15$ in 
Algorithm~\ref{alg:rominv}.

We observe that the results of the ROM based inversion are better than those of MFWI. Both the shape 
and contrasts of the inclusions are recovered well. The cross-talk from $c$ to $\rho$ is minimal. The cross-talk from $\rho$ to $c$ 
is more visible, but  small. MFWI does a worse job  in 
recovering the inclusions and minimizing the cross-talk. In particular, the cross-talk from $\rho$ to $c$ is almost as strong as the recovered $c$ inclusion.

\subsection{Marmousi model}
The objective of this experiment is to assess how the ROM inversion mitigates cycle skipping, 
which is the main cause of the spurious local minima of the conventional (MWFI) objective function.
We consider a section of the two-coefficient Marmousi model 
\cite{martin2006marmousi2}, in the domain $\Omega = [0, 6\mbox{km}] \times [0, 3.75\mbox{km}]$.
An array of  $n_s=30$ sensors is placed at the depth of $150\mbox{m}$. No noise is added to the data. The 
inversion is carried out in six layers i.e., $\ell = 6$ in Algorithm~\ref{alg:rominv}.
The optimization search space has dimension  $2 N = 2(50 \times 40) = 4000$, 
with velocity and density parameterized as in equation \eqref{eqn:rparam}.

\begin{figure}[t]
\begin{center}
\begin{tabular}{cc}
$c\; (km/s)$ & $\rho\; (g/cm^3)$ \\
\includegraphics[width=0.48\textwidth]
{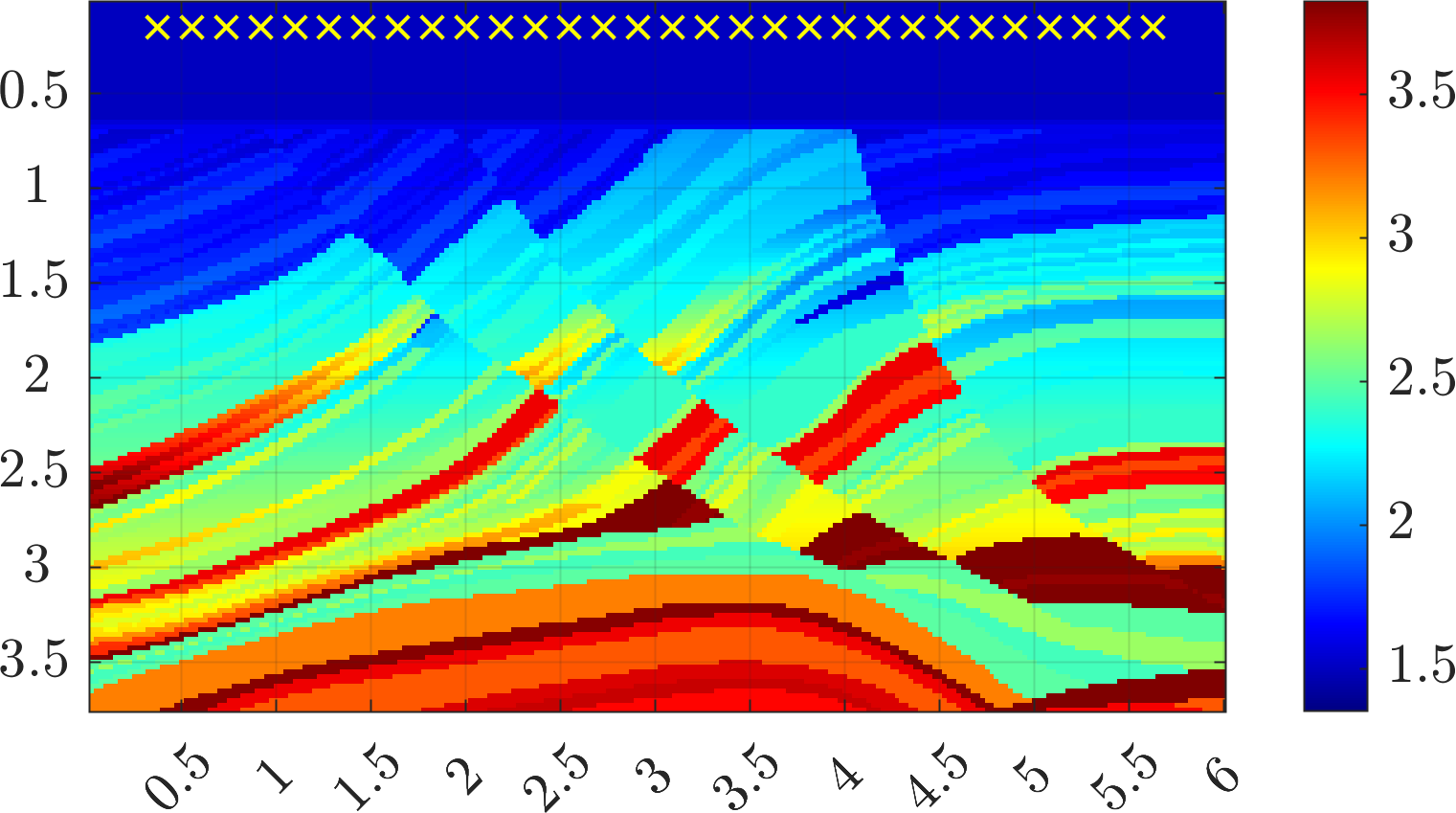} &
\includegraphics[width=0.48\textwidth]
{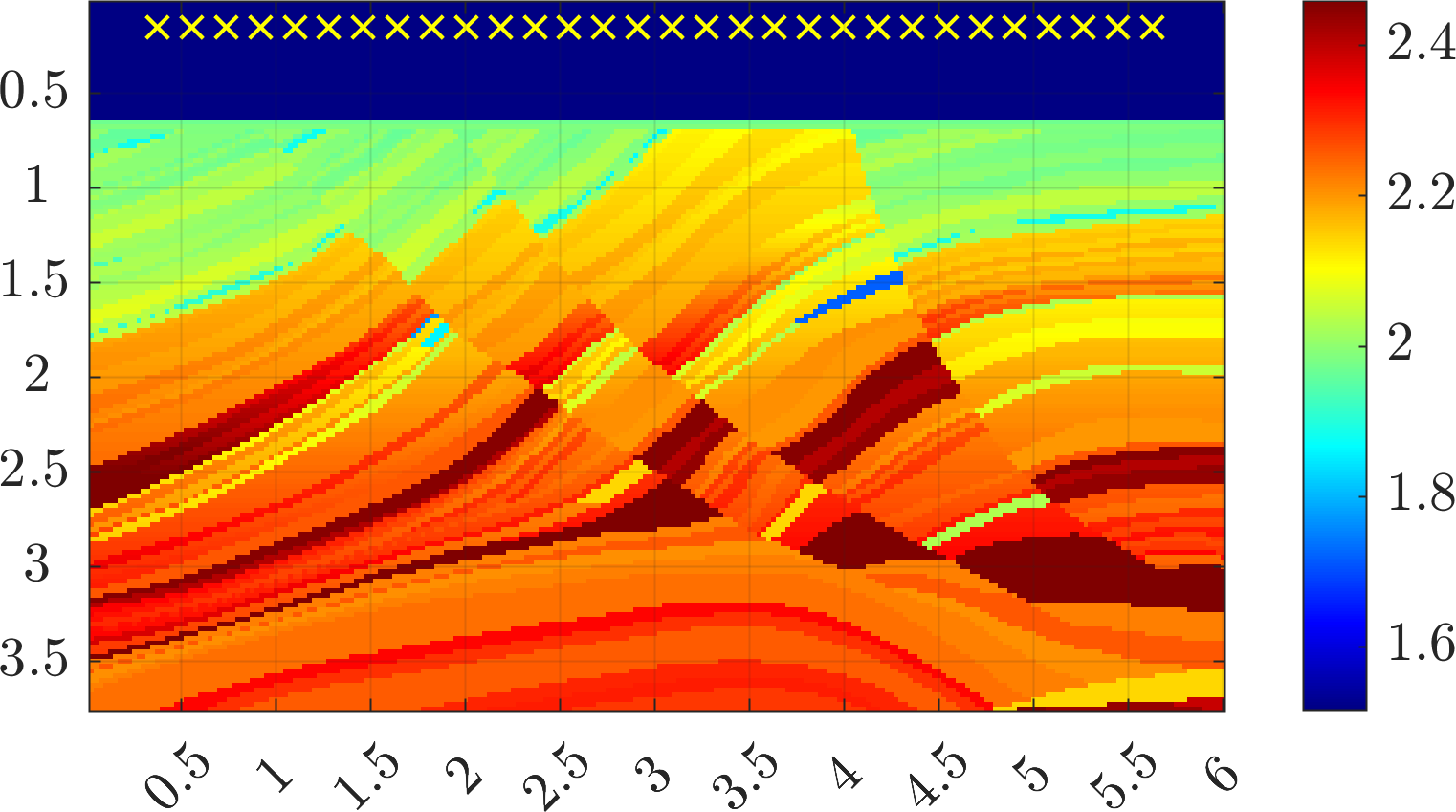} \\
\includegraphics[width=0.48\textwidth]
{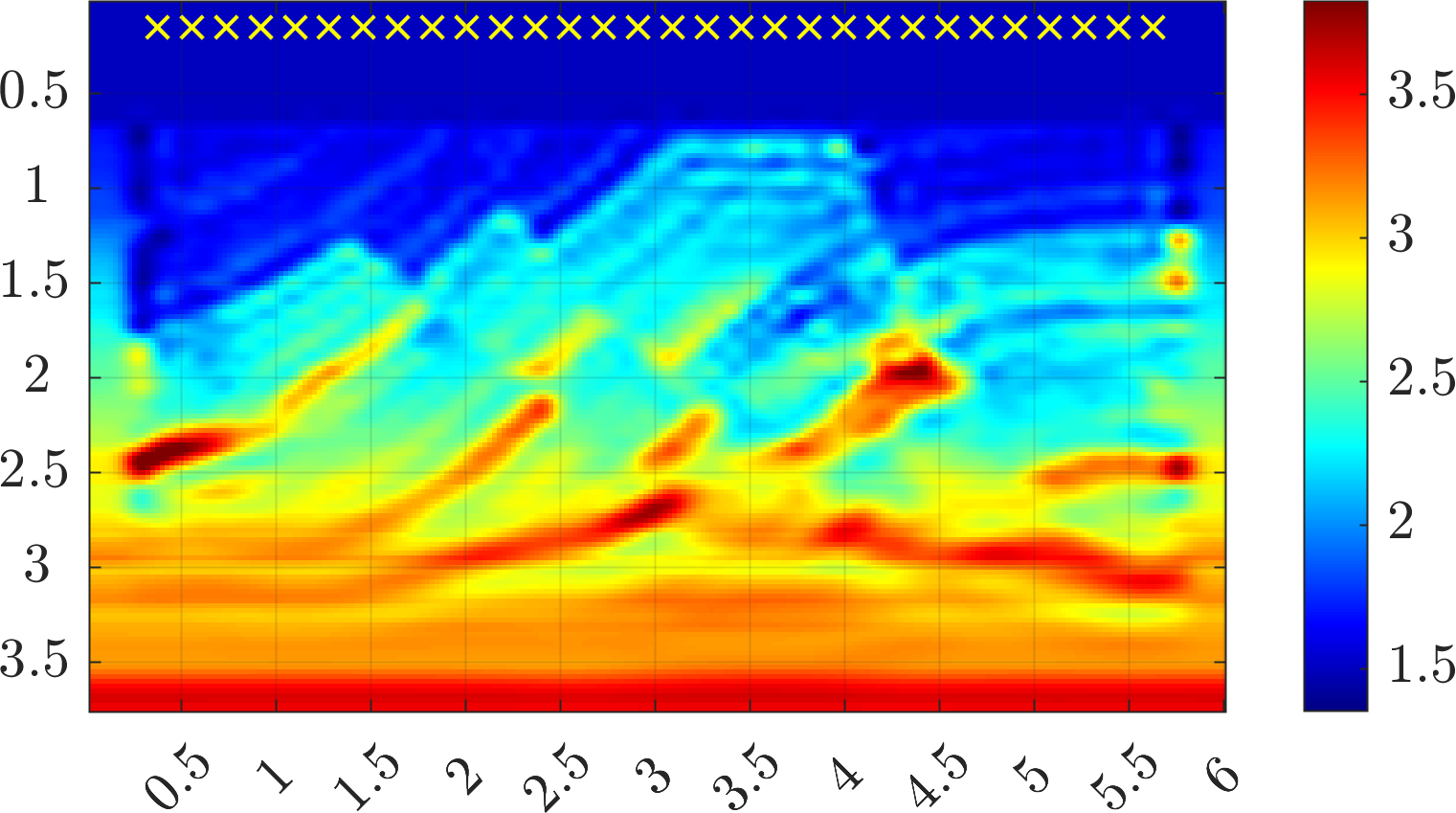} &
\includegraphics[width=0.48\textwidth]
{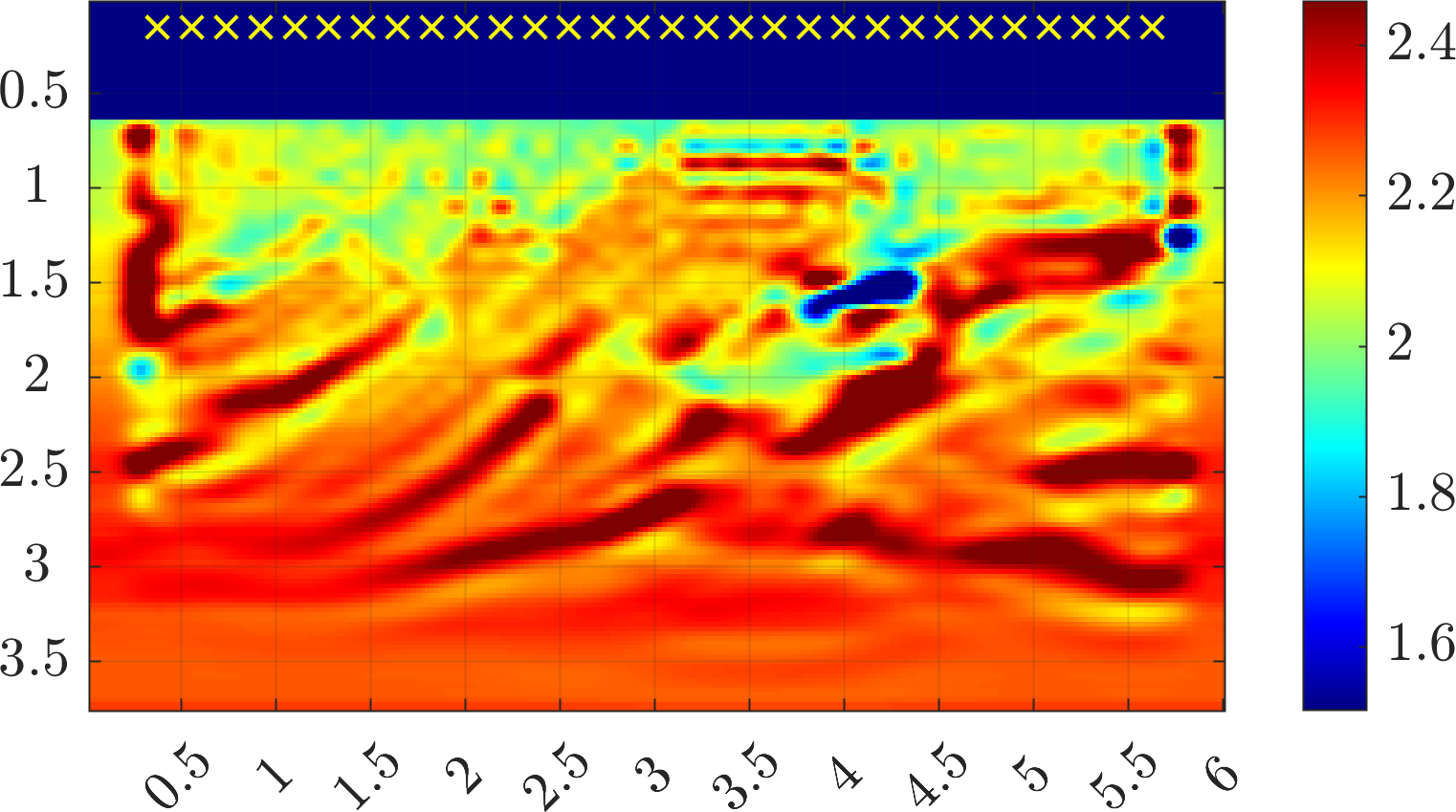} \\
\includegraphics[width=0.48\textwidth]
{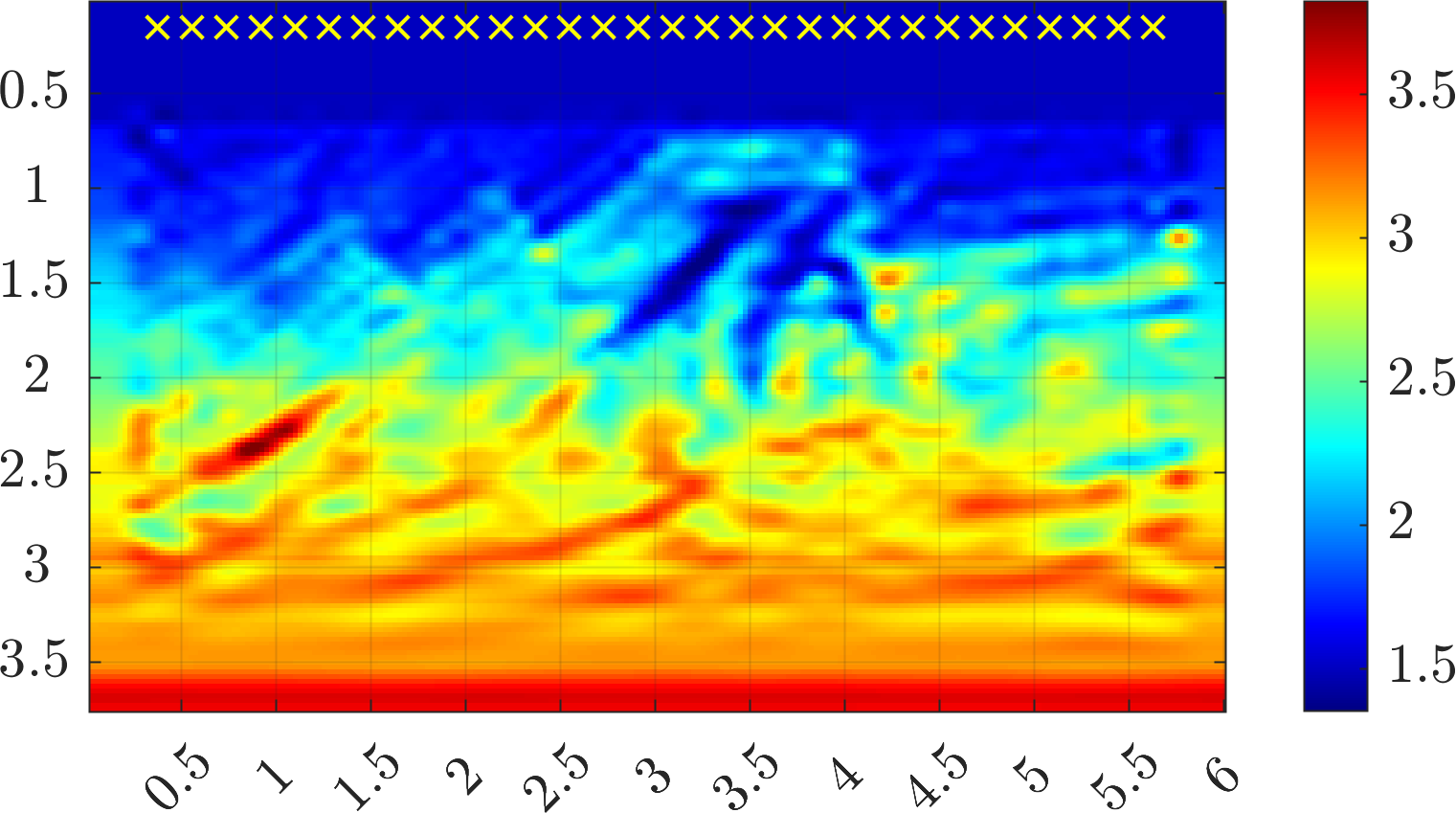} &
\includegraphics[width=0.48\textwidth]
{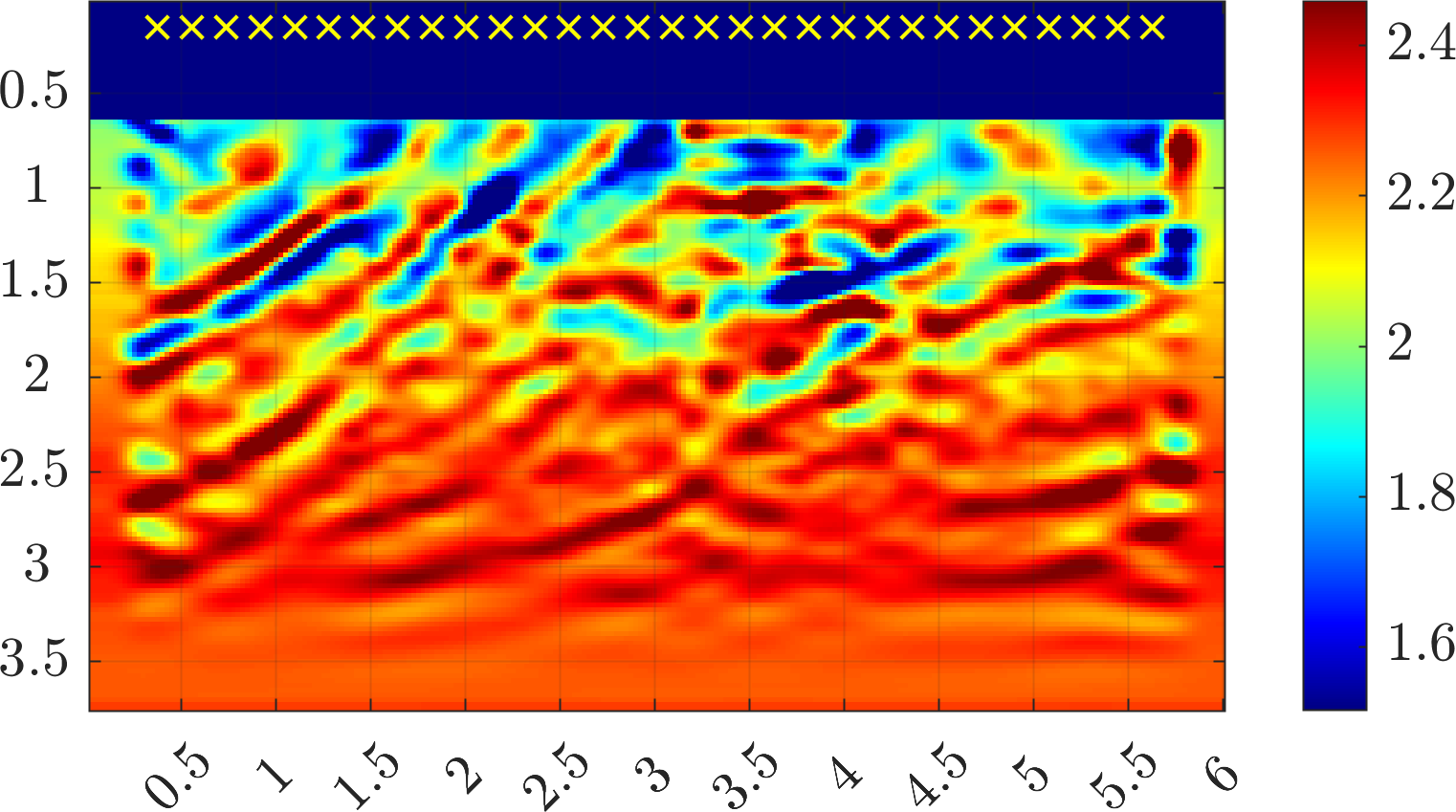}
\end{tabular}
\end{center}
\caption{Two-coefficient Marmousi model. Top row: true velocity and density; 
middle row: ROM based estimates; bottom row: MFWI estimates.
Source locations are yellow crosses. Axes are in km. The colorbar shows the contrast.
}
\label{fig:marm}
\end{figure}

The results of the second numerical experiment are shown in Figure~\ref{fig:marm}. We display both
MFWI and ROM based estimates obtained after $50$ iterations. We took $8$ iterations for the first five layers and 
10 iterations in the last layer. 
As expected, 
MFWI gets stuck in a local minimum and the result has multiple artifacts in the estimated 
velocity and a mostly nonsensical density estimate. In contrast, the ROM based estimates of $c$ and $\rho$ 
capture the true medium, except for some artifacts in the density close to the (reflective) sides of the domain. 

\subsection{Models with cracks}

In the third experiment we test the effect of noise on the inversion. See Appendix \ref{ap:A} for the noise model, with standard deviation quantified by the dimensionless parameter $b$. 

We consider two models with crack like features, embedded into smooth backgrounds.
Both models occupy the domain $\Omega = [0, 2.5\mbox{km}] \times [0, 1.2\mbox{km}]$ with variable background 
velocity ranging from $1\mbox{km/s}$ to $1.5\mbox{km/s}$ and background density in the $1-1.2\mbox{g/cm}^3$ range.
The velocity inside the cracks is $2\mbox{km/s}$ while the density is $1.5\mbox{g/cm}^3$.  The difference between the two models is that the crack features of $c$ and $\rho$ are 
interchanged.

\begin{figure}[t]
\begin{center}
\begin{tabular}{cc}
$c\; (km/s)$ & $\rho\; (g/cm^3)$ \\
\includegraphics[width=0.48\textwidth]
{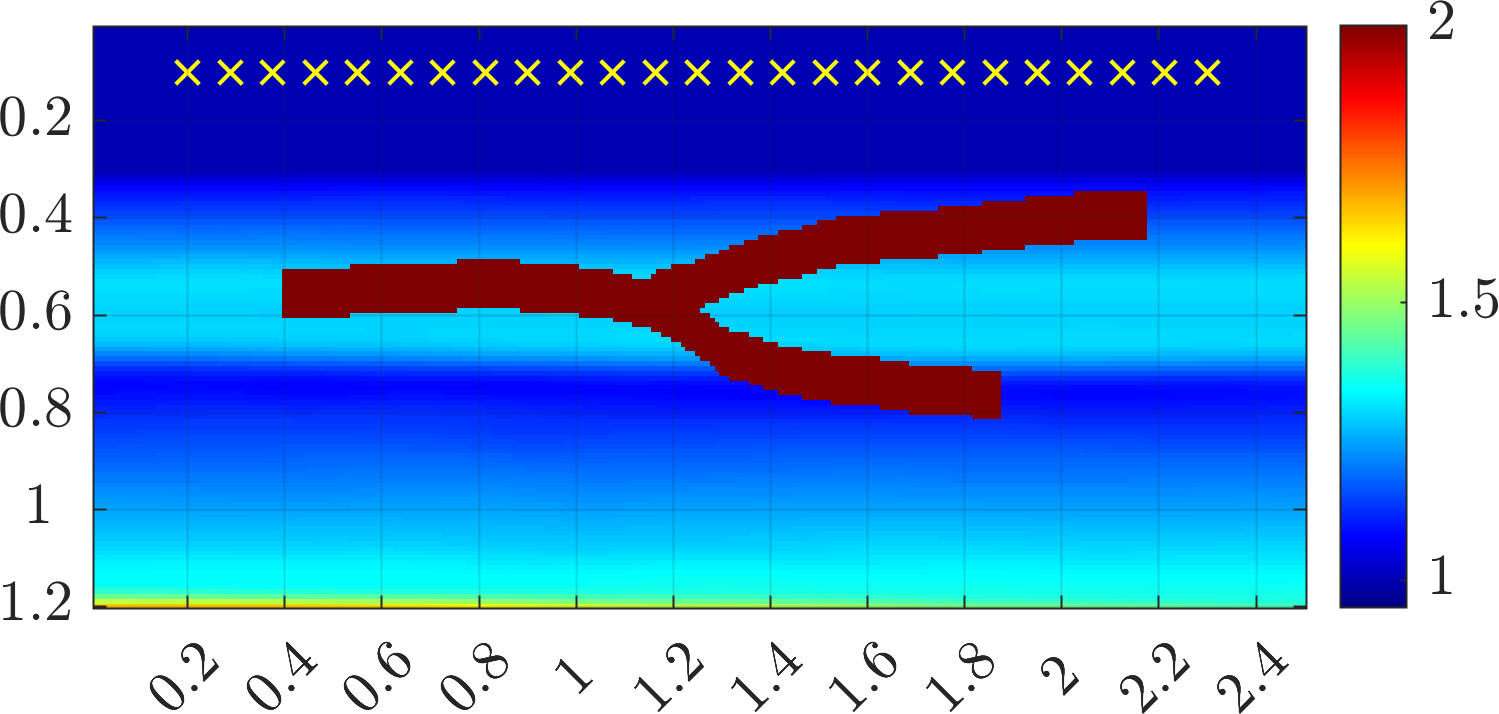} &
\includegraphics[width=0.48\textwidth]
{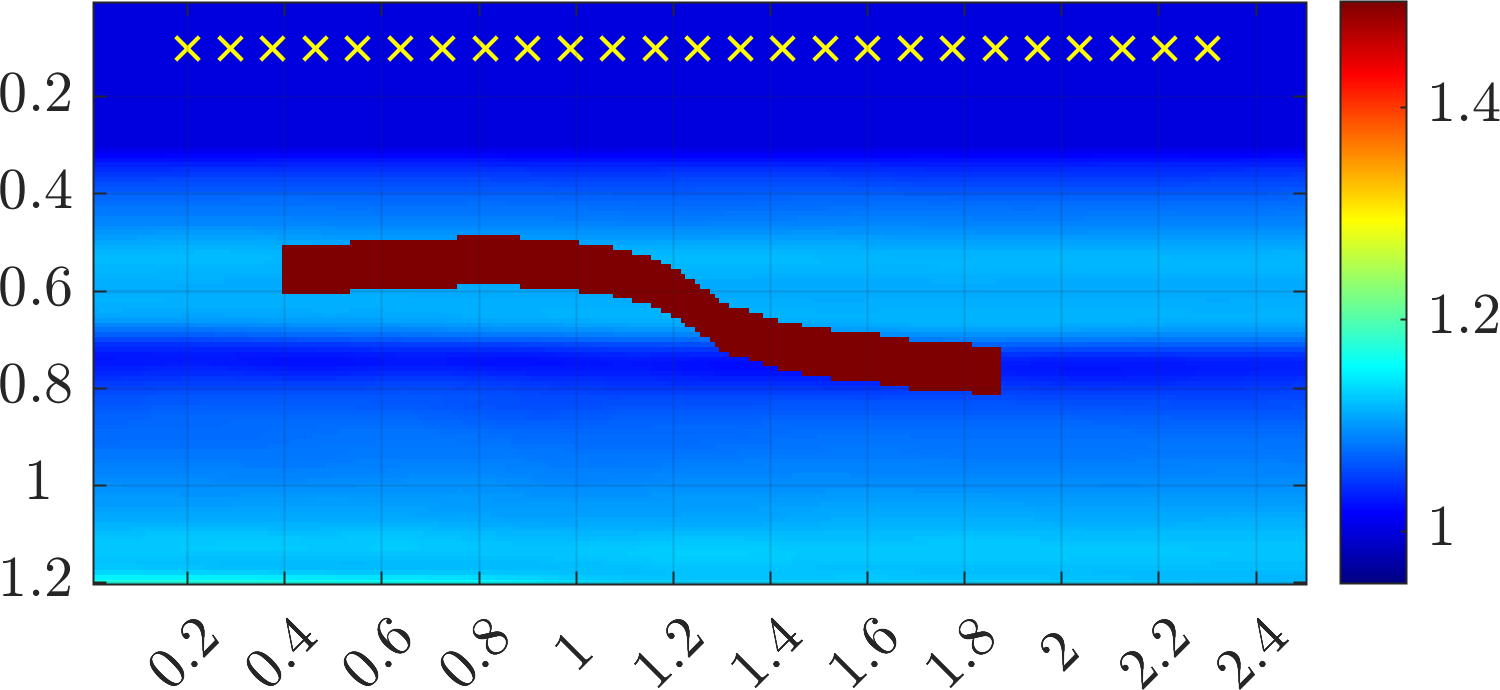} \\
\includegraphics[width=0.48\textwidth]
{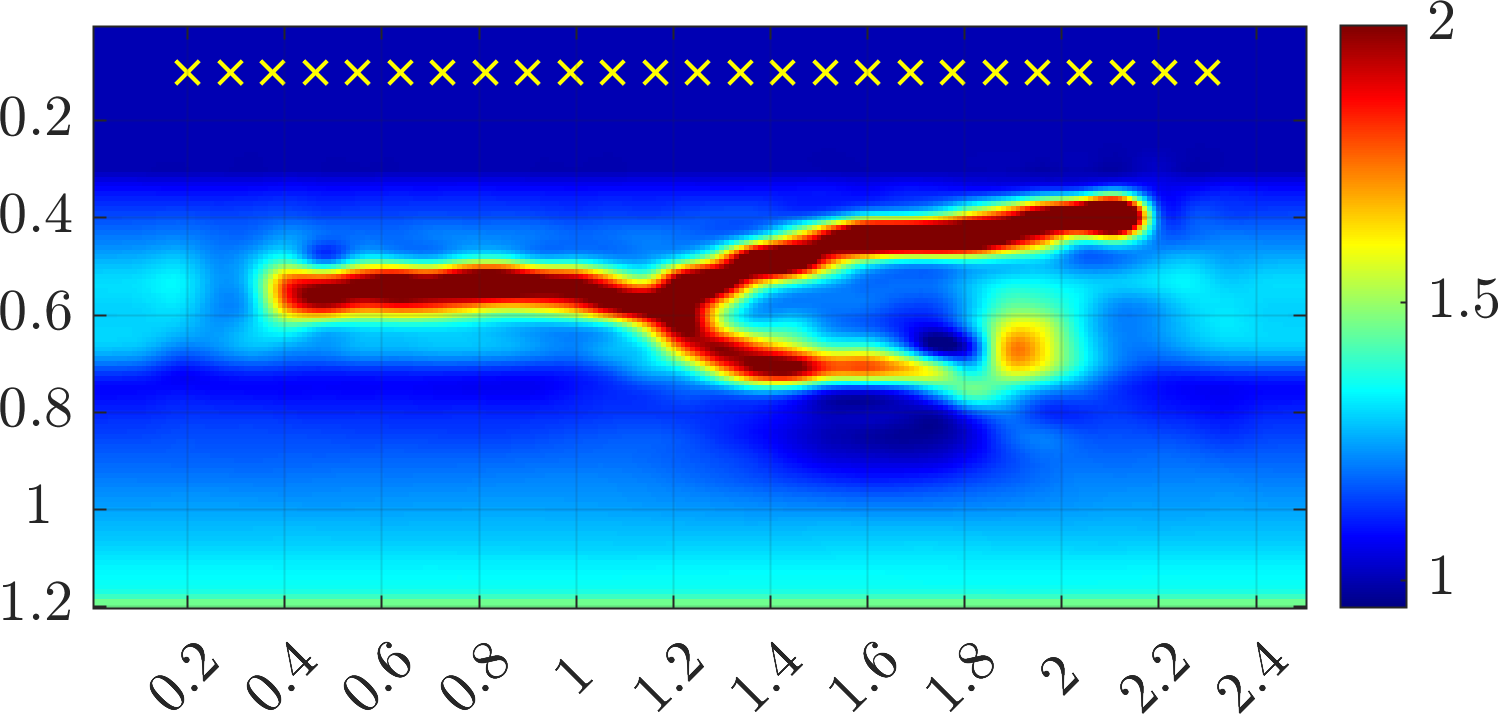} &
\includegraphics[width=0.48\textwidth]
{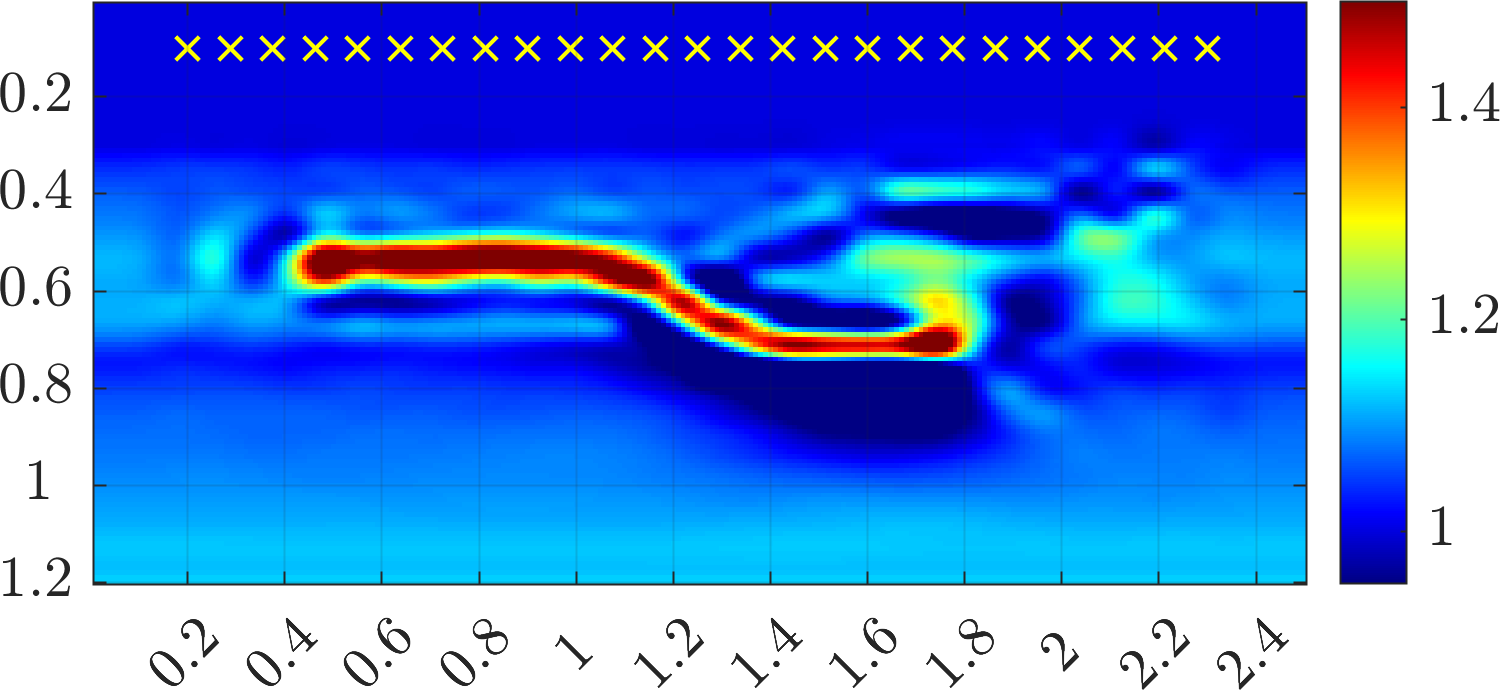} \\
\includegraphics[width=0.48\textwidth]
{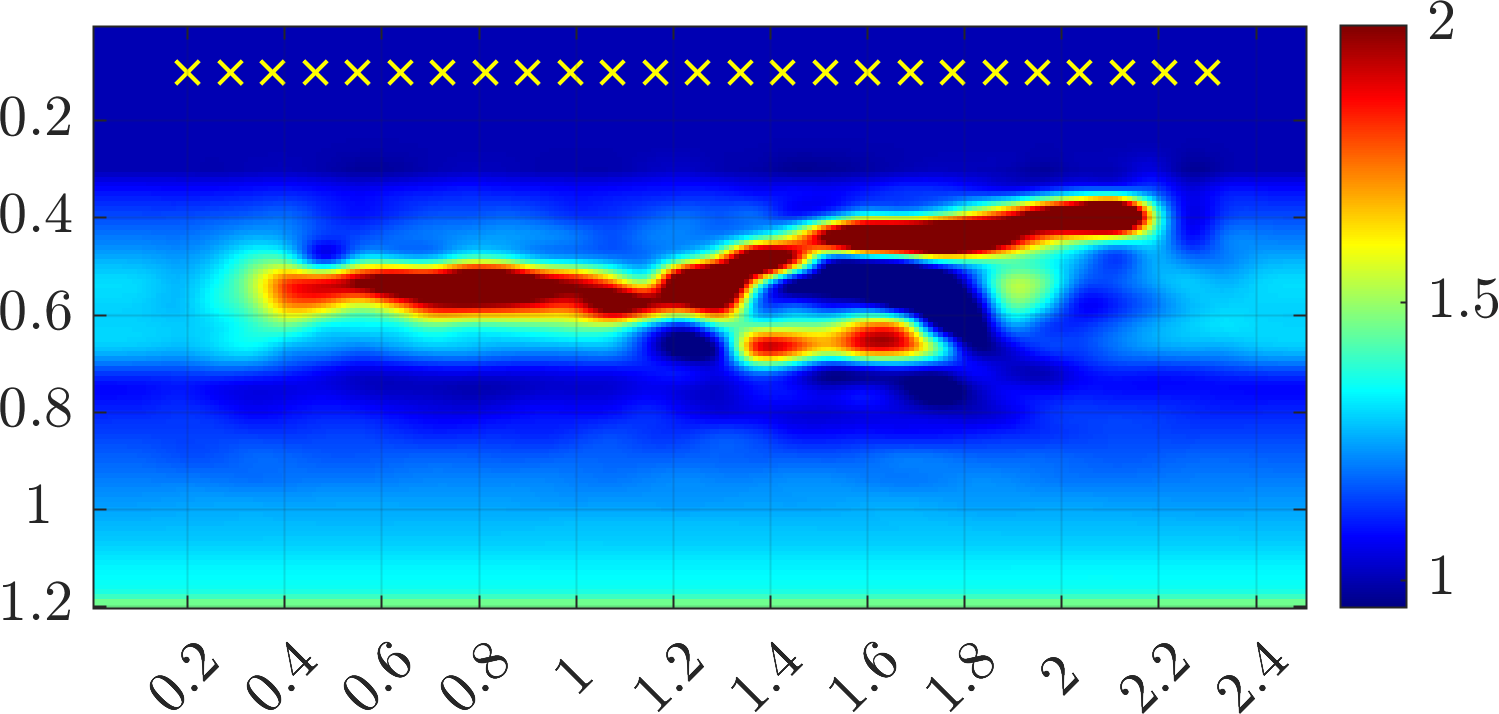} &
\includegraphics[width=0.48\textwidth]
{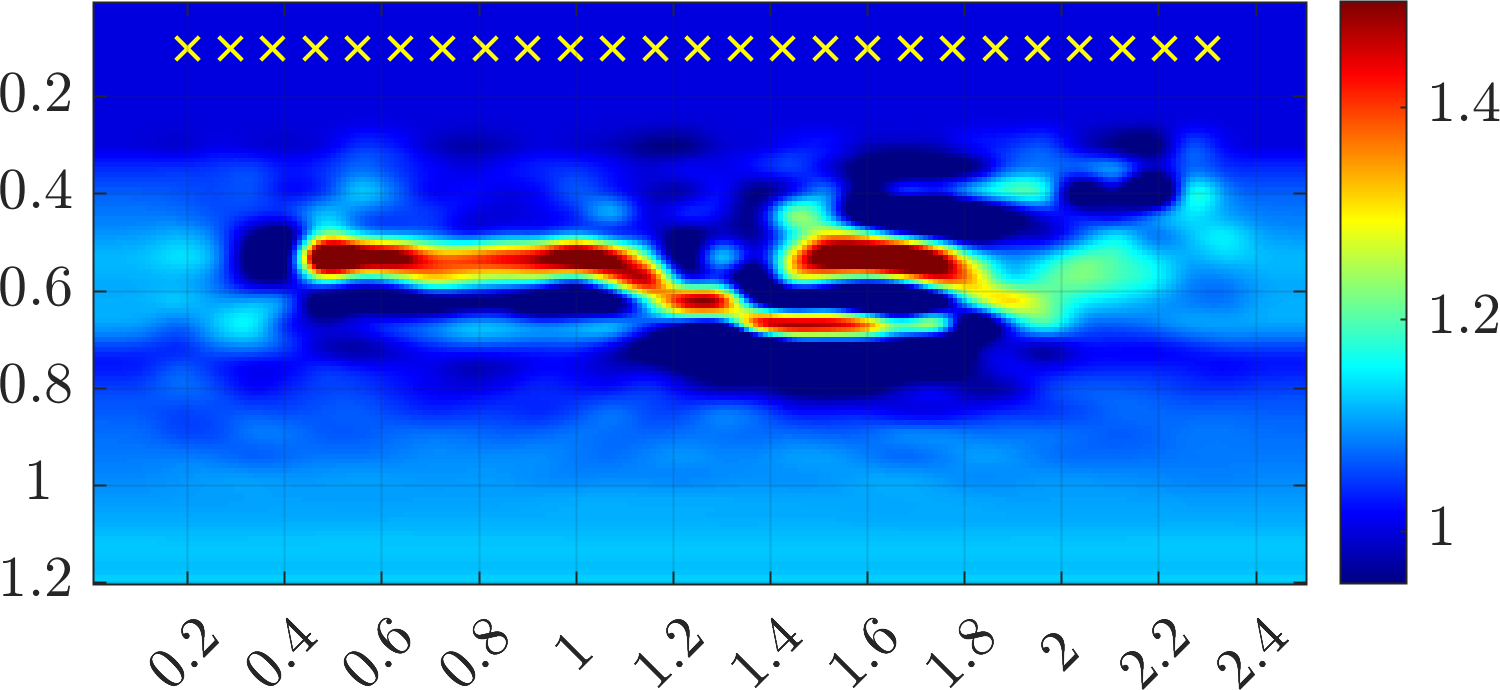} \\
\includegraphics[width=0.48\textwidth]
{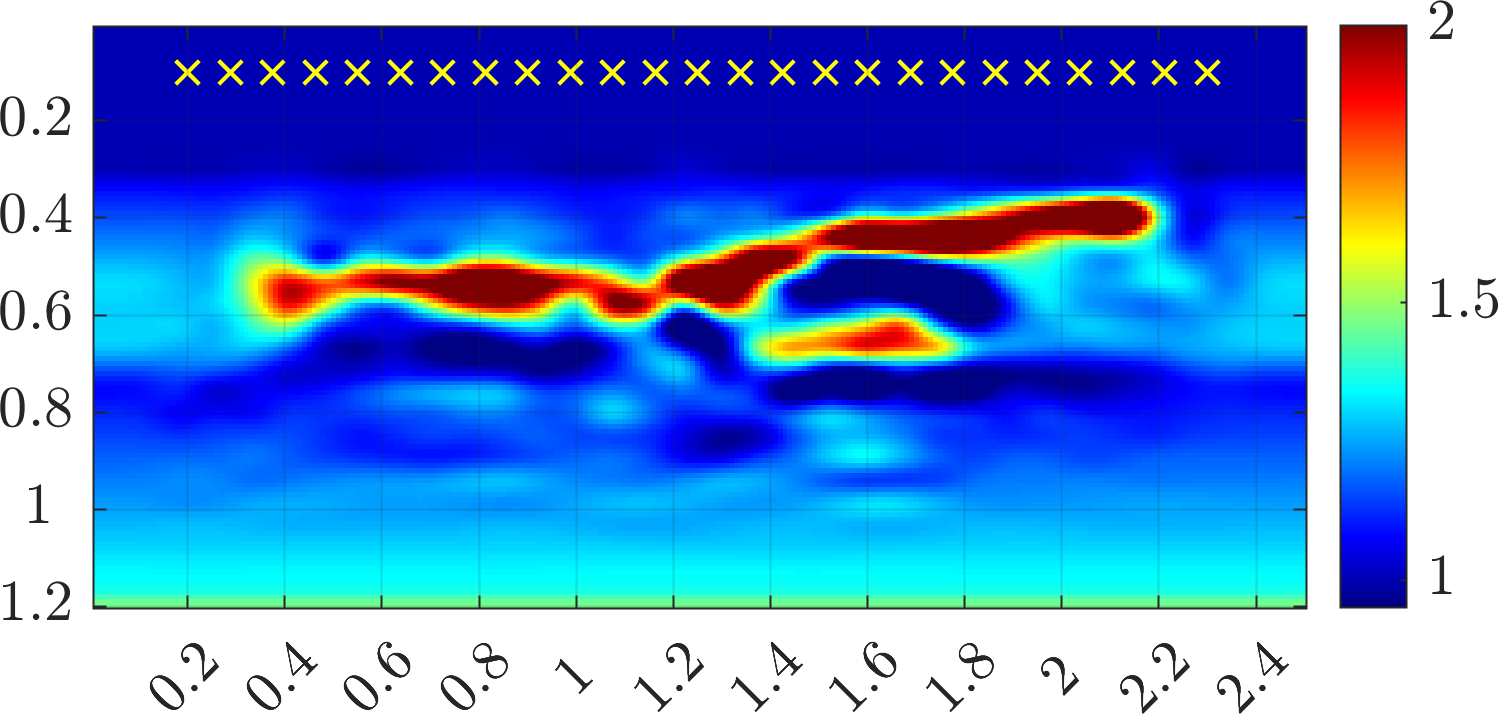} &
\includegraphics[width=0.48\textwidth]
{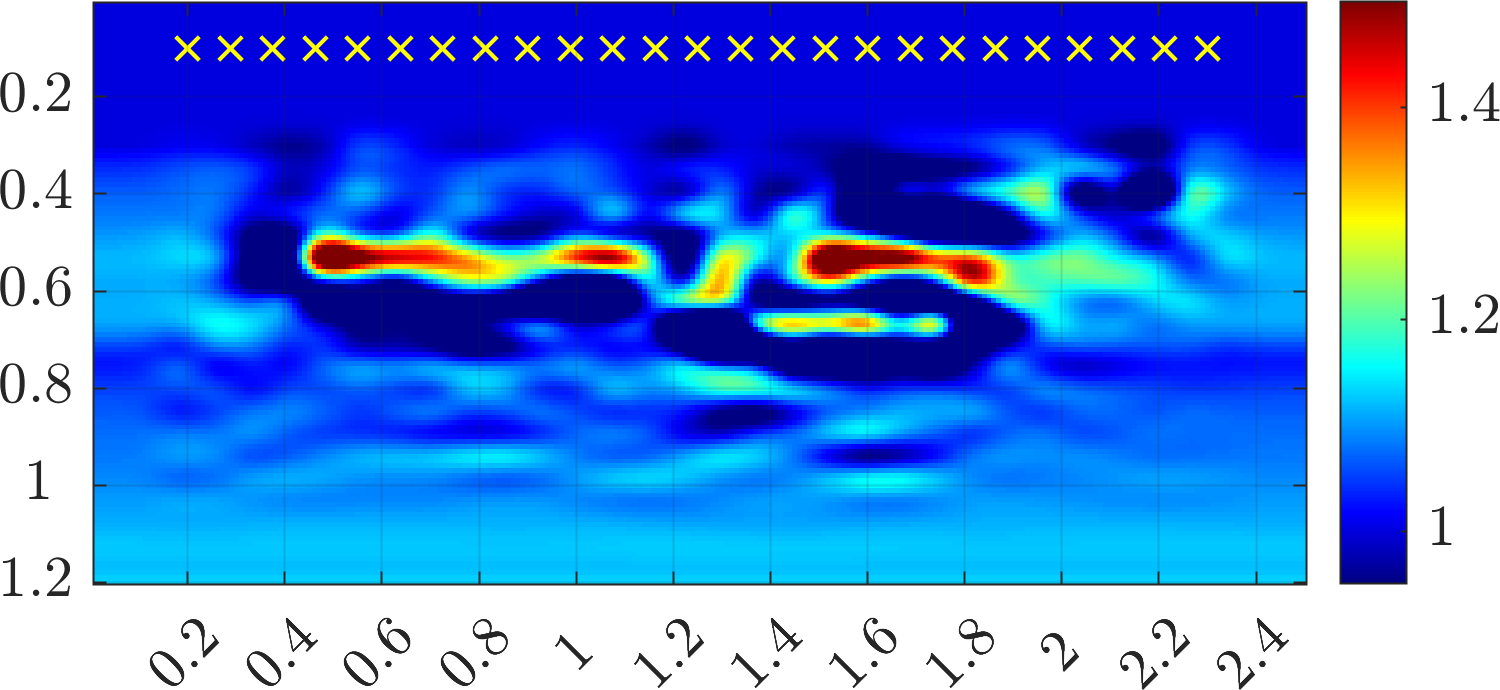} 
\end{tabular}
\end{center}
\caption{First crack model. Top row: true velocity and density; 
second row: ROM based estimates from noiseless data; 
third row: ROM based estimates from noisy data with $b = 10^{-2}$; 
bottom row: MFWI estimates from noiseless data.
Source locations are yellow crosses. Axes are in km. The colorbar shows the contrast.}
\label{fig:crack1}
\end{figure}

We invert for both the background and the cracks i.e., only the reference values  
$c_o$ and $\rho_o$ of the medium near the array are assumed known. The array of $n_s=25$ sensors is  placed at the depth of $100\mbox{m}$. 
The inversion is carried out with both noiseless and noisy data, with noise level $b = 10^{-2}$ and 
$b = 3 \cdot 10^{-2}$ for the first and second model, respectively. The optimization search space 
 has dimension $2 N = 2(30 \times 20) = 1200$, with velocity and density parameterized 
as in equation \eqref{eqn:rparam}. For both models we perform $40$ iterations of 
MFWI and the ROM based inversion.
We do not display MFWI results with noise, 
because they are not so different from the noiseless ones (MFWI is known to be robust with respect to additive Gaussian noise). The  results of the ROM based inversion with noisy data 
are shown after  $30$ iterations. A single layer ($\ell = 1$) is used in all the cases.

\begin{figure}[t]
\begin{center}
\begin{tabular}{cc}
$c\; (km/s)$ & $\rho\; (g/cm^3)$ \\
\includegraphics[width=0.48\textwidth]
{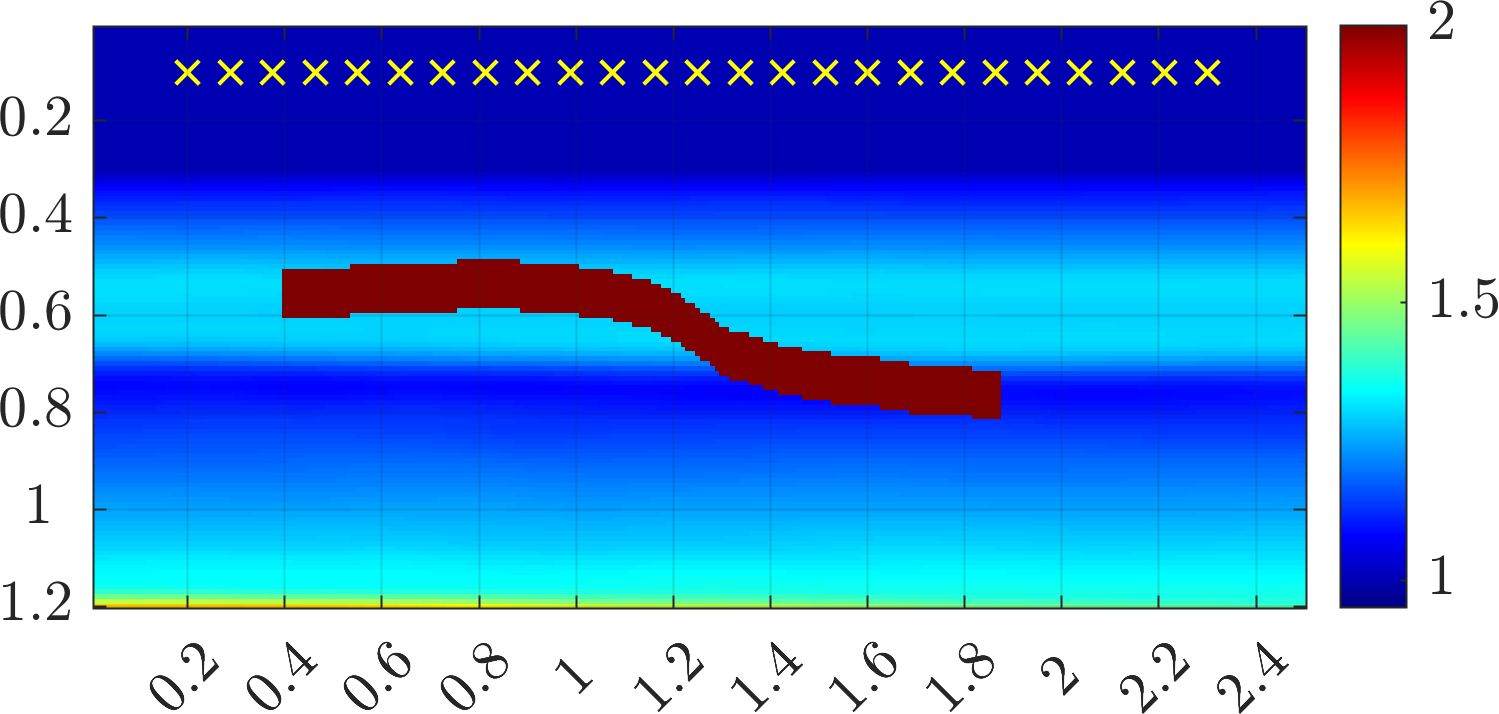} &
\includegraphics[width=0.48\textwidth]
{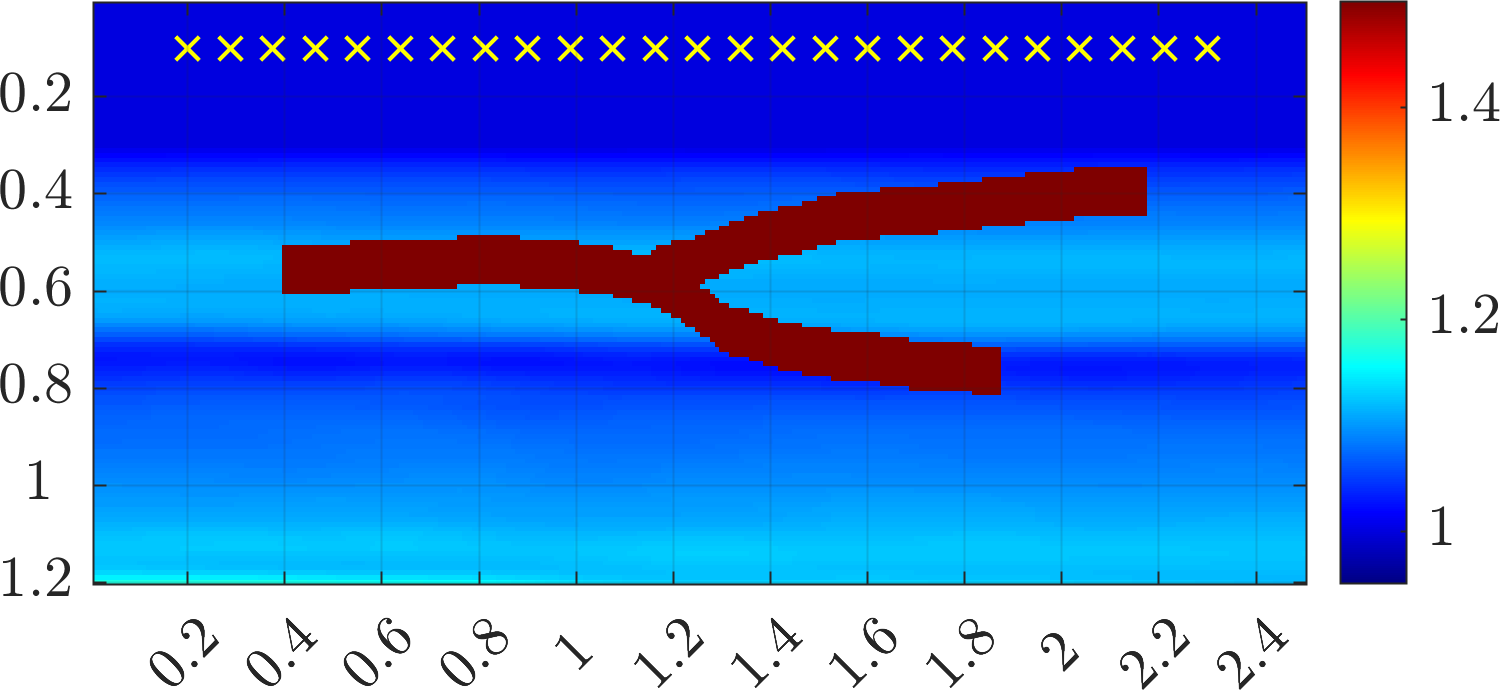} \\
\includegraphics[width=0.48\textwidth]
{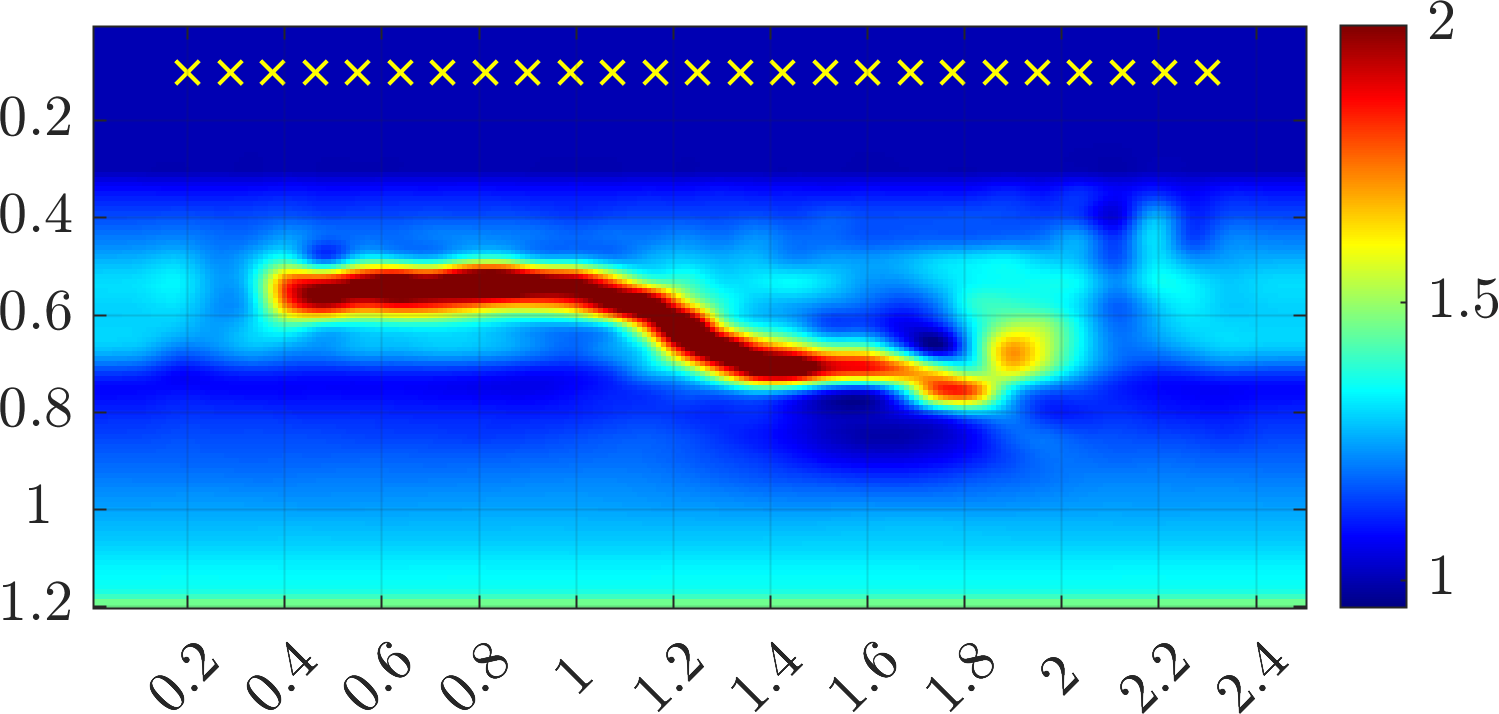} &
\includegraphics[width=0.48\textwidth]
{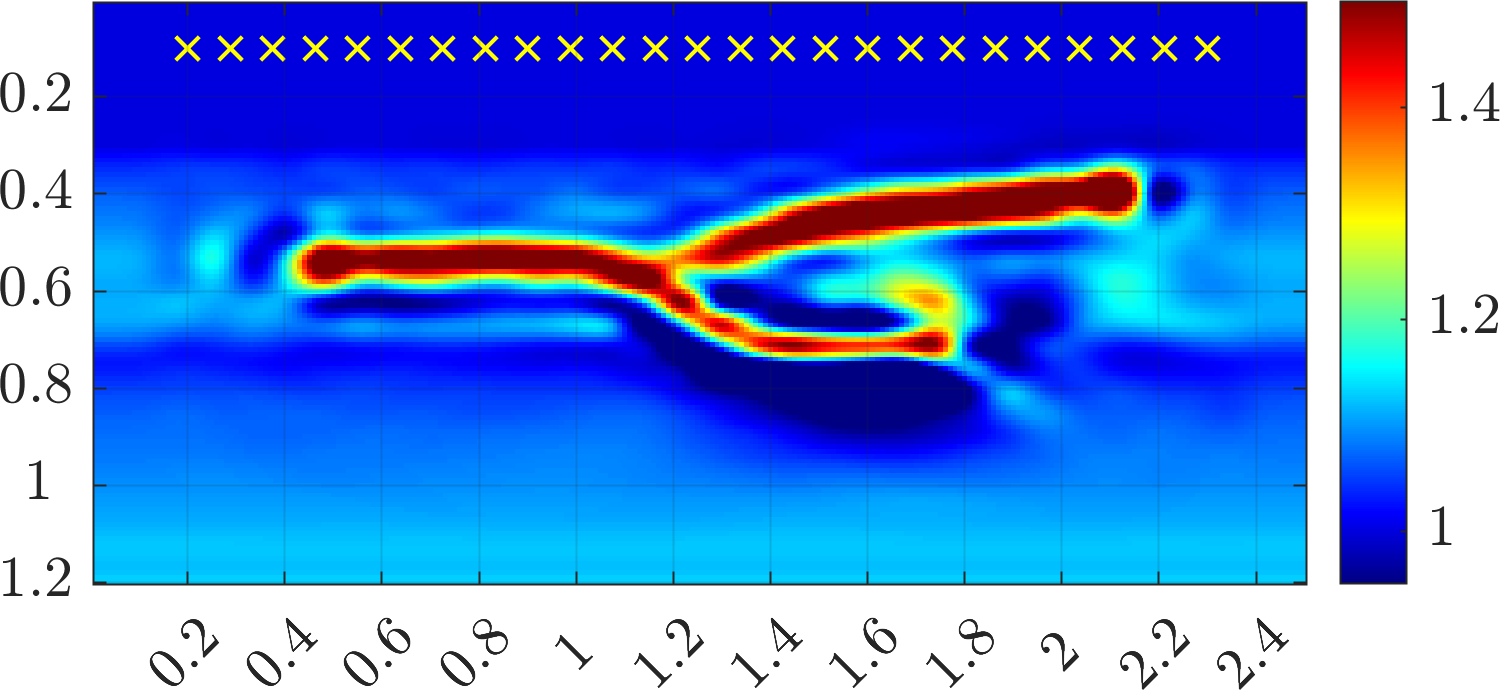} \\
\includegraphics[width=0.48\textwidth]
{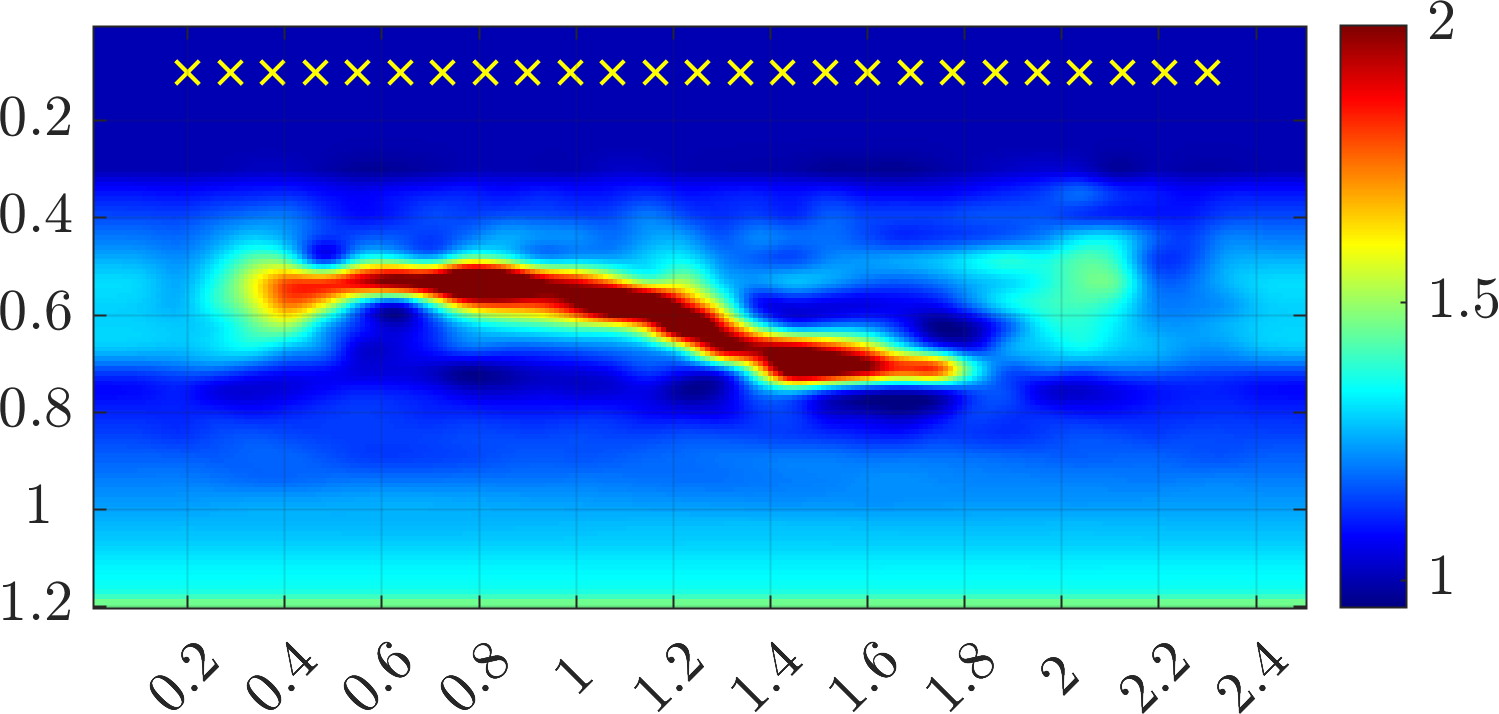} &
\includegraphics[width=0.48\textwidth]
{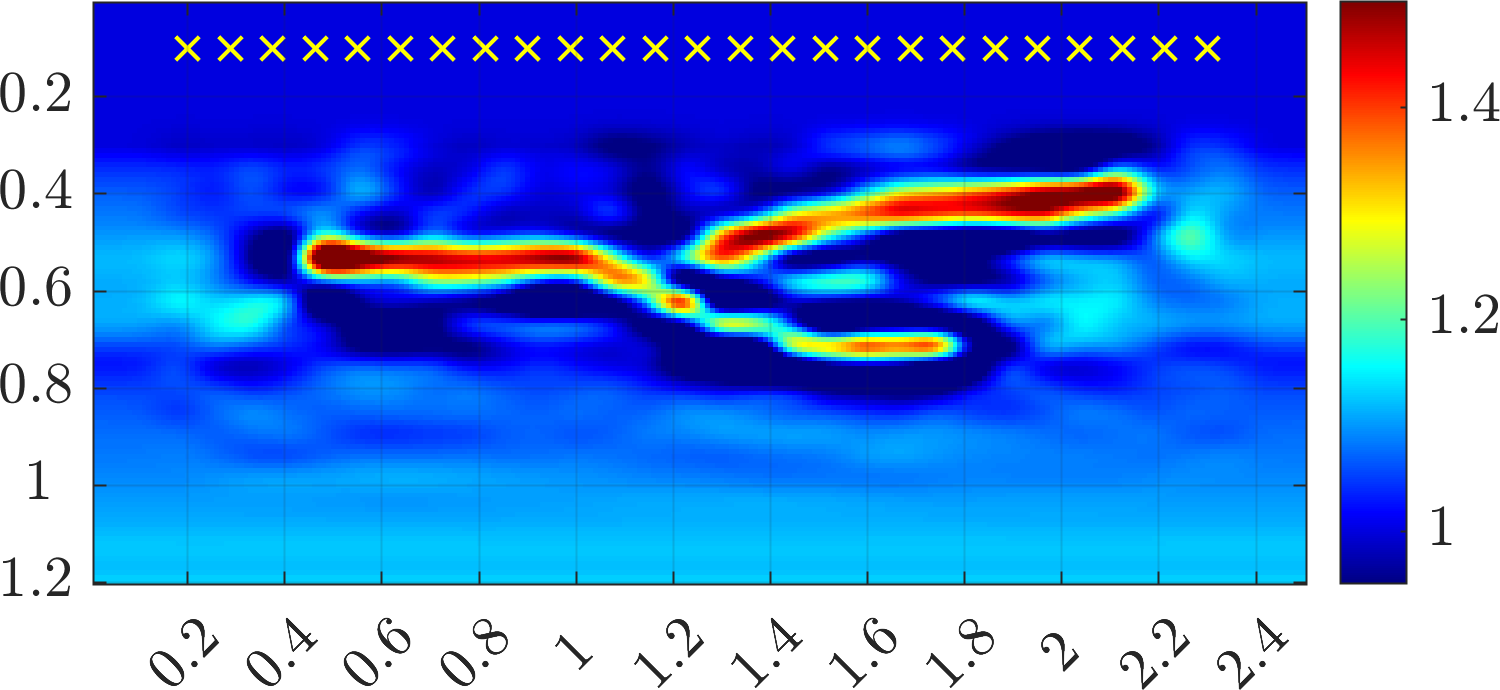} \\
\includegraphics[width=0.48\textwidth]
{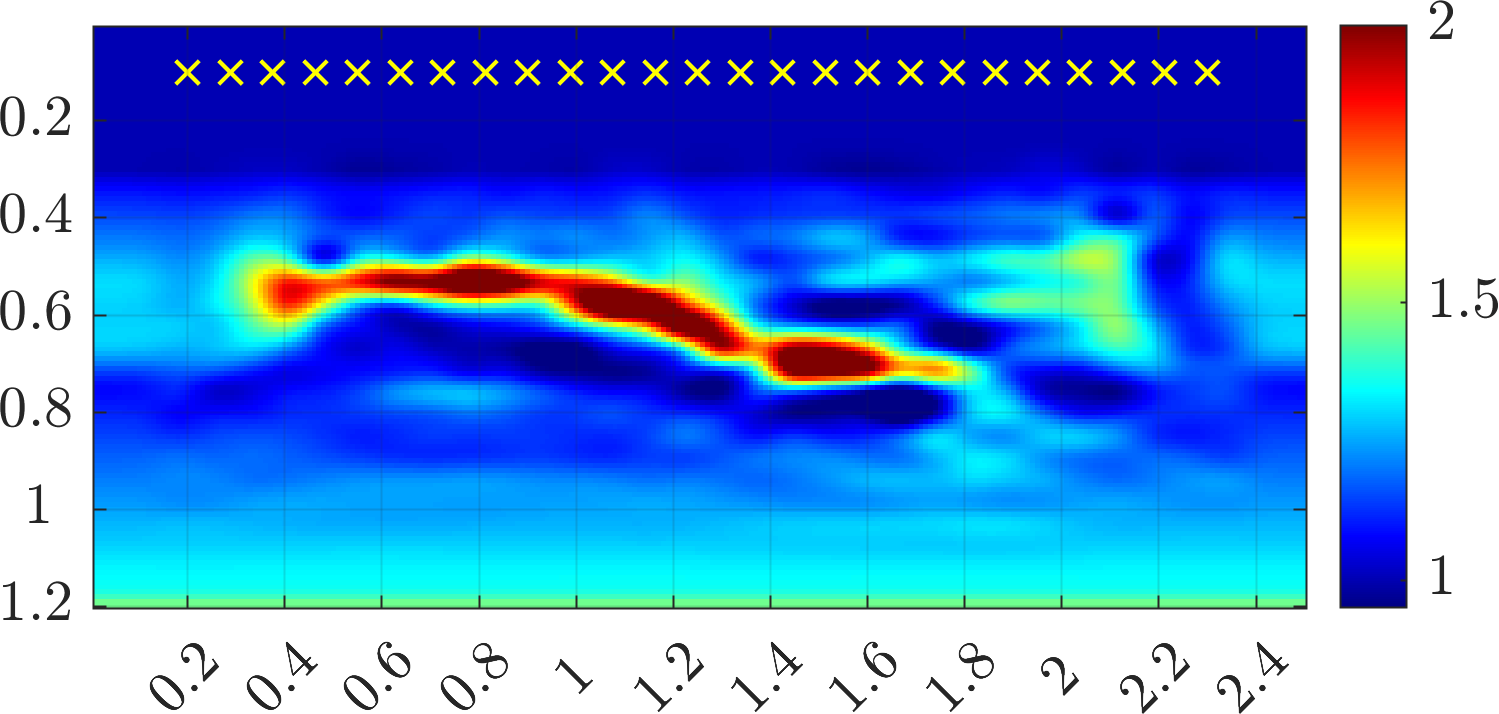} &
\includegraphics[width=0.48\textwidth]
{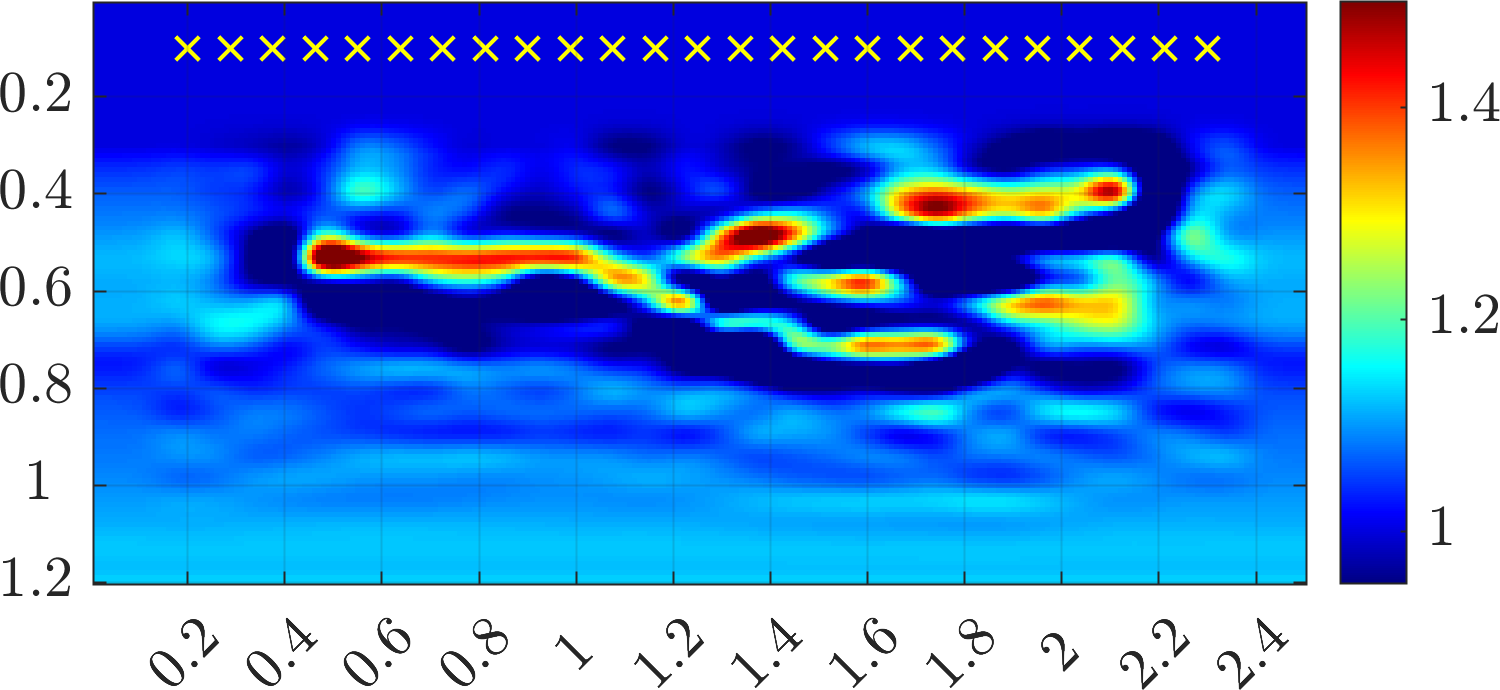}
\end{tabular}
\end{center}
\caption{Second crack model. Top row: true velocity and density; 
second row: ROM based estimates from noiseless data; 
third row: ROM based estimates from noisy data with $b = 3 \cdot 10^{-2}$; 
bottom row: MFWI estimates from noiseless data.
Source locations are yellow crosses. Axes are in km. The colorbar shows the contrast.}
\label{fig:crack2}
\end{figure}

The results for the first model are  in Figure~\ref{fig:crack1}. 
The bottom branch of the crack feature in $c$ is difficult to recover, because most of the energy 
is reflected back to the array by the top branch. 
MFWI gives a worse recovery of this branch and also introduces an artifact in the density estimate. The ROM based estimate from noisy data is 
qualitatively similar to that of MFWI, which implies that the information about the deeper part of the 
model is masked by the noise and is discarded during the regularization process.

The results of the experiments for the second model with cracks are shown in Figure~\ref{fig:crack2}.
This model is easier to recover due to the fork-like crack being in the density part of the model rather 
than in the velocity part. Thus, we observe that the cracks are estimated well  by the ROM based 
inversion both with noisy and noiseless data. MFWI gives a worse estimate of 
the density crack even from the noiseless data.
\section{Summary}
\label{sect:Summary}
We introduced a data driven reduced order model (ROM) for a general first order hyperbolic system of equations 
that governs the evolution of all linear types of waves (acoustic, electromagnetic and elastic) in a lossless medium occupying some domain $\Omega$. 
We derived an exact time stepping scheme for the hyperbolic system, on a uniform time grid. This scheme maps the wave field snapshots from one instant to the next, using a unitary operator called the wave propagator.
The ROM is defined as an orthogonal, Galerkin  projection of the time stepping scheme on the space spanned by the wave snapshots.
We are interested in a ROM that can be used to solve inverse wave scattering problems, where the snapshots are unknown
almost everywhere in $\Omega$. The ROM is called data driven, because it can be computed from inverse scattering data i.e., measurements of the snapshots at 
a few sensor locations which probe the medium in $\Omega$ with short signals (pulses) and record the generated, backscattered waves. The mapping between these  measurements  and the ROM is nonlinear, but we show that it can be computed using efficient techniques from numerical linear algebra. 

The definition of the ROM respects the causality of wave propagation: The matrix of ROM snapshots is block upper triangular and the ROM propagator, the 
Galerkin projection of the propagator operator,  is a block upper Hessenberg matrix.  
The size of the blocks of these matrices is given by the number of wave excitations. The block algebraic structure reflects the causal dependence on the data: The  top left $k$ blocks are determined by the measurements up to the $k^{\rm th}$ time instant. This is essential not only because the ROM captures correctly
the physics of wave propagation, but also because we can use it to solve the inverse problem in a layer peeling fashion. 

While the definition of the ROM is general, we specialize its application to inverse scattering with acoustic waves. The wave field consists 
of the acoustic pressure and velocity and the medium in $\Omega$ is modeled by the unknown and variable wave speed $c$
and density $\rho$. We use the ROM to approximate the mapping from the data to the vectorial wave field at inaccessible points in $\Omega$, aka the internal wave. 
This wave is computed for the best guess of $c$ and $\rho$,  and it  fits the data by design. However,  it is not a solution of the hyperbolic system of equations unless the guess of the medium is right. This motivates the formulation of  the inverse problem 
as a minimization of  the discrepancy between the approximated wave field and the solution of the hyperbolic system of equations.
Carefully designed numerical simulations show that this formulation consistently outperforms the minimization of the nonlinear least squares data misfit objective function, used by the existing inverse scattering methodology.

\bc{The ROM construction and inversion methodology is specialized to the time domain formulation of wave propagation. There are ROMs for frequency domain formulations, with array measurements of time harmonic waves at multiple frequencies. Examples  of such ROMs are in \cite{borcea2014model}  for parabolic equations and in \cite{borcea2021reduced} for waves  in lossy one dimensional media. The construction of these ROMs is different and they do not enjoy the useful causality properties of the ROMs built in the time domain.  The frequency domain formulation is important for dealing with dispersive and lossy media, which motivates future work. }

\section*{Acknowledgements}
This material is based upon research supported in part by the ONR award number N00014-21-1-2370, the AFOSR award number FA9550-22-1-0077, the NSF award number 2309197 and  by Agence de l’Innovation de D\'efense—AID via Centre Interdisciplinaire d’Etudes pour la D\'efense et la S\'ecurit\'e—CIEDS—(project PRODIPO).

\appendix
 \section{Setup for the numerical simulations}
 \label{ap:A}%
The numerical simulations are carried out 
in a two-dimensional ($d = 2$) rectangular domain $\Omega$, with sound soft boundary condition \eqref{eq:Ac6}. We use $n_{\cE} = 2 n_s$ excitations $\bF_{\be}(\bx)$ of the form \eqref{eq:Ac4}, with sources located at 
$\bx_{\eps_1}$, $\eps_1 = 1,\ldots,n_s$, and with two polarizations $\eps_2 \in \{ 1, 2 \}$ at each location. 
The probing pulse is
\begin{equation}
s(t) = 1_{[-t_s,t_s]}(t) \frac{d}{dt} \left( \cos (2 \pi \nu t) \exp \left[ - \frac{(2 \pi B)^2 t^2}{2} \right] \right),
\end{equation}
with central frequency $\nu = 6$Hz and bandwidth $B = 4$Hz. Here $1_{[-t_s,t_s]}(t)$ 
is the indicator function of the interval $[-t_s,t_s]$, where $t_s = 1.5 / \big( \nu + B \big) = 150$ms.  
The noiseless data \eqref{eq:LS12} are generated by solving numerically \eqref{eq:LS9} using an 
exponential integrator\footnote{ The derivation of the ROM uses that the wave field is a time-dependent flow, as stated in equation \eqref{eq:LS10}. This is why we used a forward solver 
 that preserves this property. } based on \cite{al2011computing} with a time step $\tau_f$ defined below. 
The operator \eqref{eq:Ac9} is discretized using a finite difference scheme with two-point first-order 
discretization of all the partial derivatives. The initial state \eqref{eq:LS8} is calculated using a backward 
Euler scheme with the same time step $\tau_f$ as in the exponential integrator. 
The data are sampled at intervals $\Delta t = 1 / \big[2.3 \big( \nu + B \big) \big]= 43.5$ms. The numerical time integration 
step is $\tau_f = \Delta t / 10$.

It remains to describe the noise model: Consider the noise signal
${\mathfrak{N}}(t) \in \RR^{n_{\cE} \times n_{\cE}}$ with $t \in [-t_s,  t_{\max} + t_s]$ discretized
on a uniform time grid with step $\tau_f$, where $t_{\max} = (n_t-1)\Delta t$. 
The entries of ${\mathfrak{N}}(t)$ are identically distributed,  Gaussian, with mean zero and standard 
deviation chosen so that 
\begin{equation}
\bc{\Big( \int_{0}^{t_{\max}} dt \left\| {\mathfrak{N}}(t) \right\|_F^2 \Big)^{1/2} = 
b \Big( \int_{0}^{t_{\max}} dt \left\| \widetilde{\bA}(t) \right\|_F^2 \Big)^{1/2},}
\label{eqn:noiselevel}
\end{equation}
for a user defined noise level $b$,
where
\begin{equation}
\bc{\widetilde{\bA}_{\be',\be}(t) = \int_{\Omega} d \bx \, \begin{pmatrix}\bF_{\be'}(\bx) \\ 0 \end{pmatrix}^T \bpsi_{\be}(t,\bx) \stackrel{\eqref{eq:Ac8}}{\approx} \sqrt{\zeta_o} {\itbf e}_{\eps'_2}^T \bu_{\be}(t,\bx_{\eps'_1}).}
\end{equation}
The noisy data $\bD_j^{\mathfrak{N}}$ are 
\begin{equation}
\bD_j^{\mathfrak{N}} = \bD_j + \int_{- t_s}^{t_s} dt \, {\mathfrak{N}}(t_j -t)s(-t), \quad
j = 0,\ldots,n_t-1,
\label{eqn:noisydata}
\end{equation}
where the integrals in \eqref{eqn:noiselevel} and \eqref{eqn:noisydata} are computed numerically using 
the trapezoidal rule on the same grid on which ${\mathfrak{N}}(t)$ is defined.

\bc{The regularization parameter $r$ can be determined with a trial-and error approach: One can start with a large $r$ and then see how the results improve or deteriorate if we decrease it. This is sufficient in the examples addressed in this paper.
A more systematic approach can be based on estimating the noise level from
the deviation of the data from the reciprocity relation. We refer to Appendix F in \cite{borcea2022waveform} for the description of such an approach.}

 \section{The block Arnoldi iteration}
 \label{ap:B}
 The following algorithm computes the orthogonal transformation that puts the regularized ROM propagator in the 
 causal, block upper Hessenberg form. We write it for a generic matrix ${\itbf X}$ with $m\times m$ blocks.

\vspace{0.05in} 
\begin{algorithm}\textbf{\emph{(block Arnoldi iteration)}}
\label{alg:blockarnoldi}%
\vspace{0.04in} \noindent \\ \textbf{Input:} The matrix ${\itbf X} \in \RR^{mn \times mn}$ and the starting block 
$\by \in \RR^{mn\times m}$.

\vspace{0.04in} \noindent 
Compute ${\itbf q}_0 = \by (\by^T \by)^{-1/2} \in \RR^{mn \times m}$;\\ 
For $k = 1,\ldots,n-1$:\\
- Set ${\itbf w} = {\itbf X} {\itbf q}_{k-1} \in \RR^{mn \times m}$;\\
- Orthogonalize 
${\itbf w} = {\itbf w} - [{\itbf q}_0,\ldots,{\itbf q}_{k-1}] \left([{\itbf q}_0,\ldots,{\itbf q}_{k-1}]^T {\itbf w}\right)$;\\
- Normalize ${\itbf q}_{k} = {\itbf w} ({\itbf w}^T {\itbf w})^{-1/2} \in \RR^{mn \times m}$.

\vspace{0.04in} \noindent  \textbf{Output:} Orthogonal matrix 
${\itbf Q} = [{\itbf q}_0,\ldots,{\itbf q}_{n-1}] \in \RR^{mn \times mn}$ such that
${\itbf Q}^T {\itbf X} {\itbf Q} \in \RR^{mn \times mn}$ is block upper Hessenberg.
\end{algorithm}\\

\vspace{0.02in} \noindent 
To obtain $({\itbf w}^T {\itbf w})^{-1/2}$, we use the unique, symmetric positive definite matrix square root computed via the spectral decomposition of the Grammian ${\itbf w}^T {\itbf w}$.

\bibliographystyle{siam} 
\bibliography{biblio.bib}

\end{document}